\def\ung{{{\frak{g}}}}
\def\ungh{{{{\hat{\frak{g}}}}}}
\def\uqg{{{U_{q}(\ung)}}}
\def\uqres{U_q^{\roman{res}}(\ung)}
\def\uqresh{U_q^{\roman{res}}(\hat{\ung})}
\def\Tau{{\Cal T}}
\def\uqresho{U_q^{{\roman{res}}}(\hat\ung)^0}
\def\uqreshm{U_q^{{\roman{res}}}(\hat\ung)^-}
\def\uqreshp{U_q^{{\roman{res}}}(\hat\ung)^+}
\def\uepresh{U_\epsilon^{\roman{res}}(\hat{\ung})}

\def\uepresho{{U}_{\epsilon}^{{\roman{res}}}(\hat\ung)^0}
\def\oh{\overline{h}}
\def\oe{\overline{e}}
\def\uepres{U_\epsilon^{\roman{res}}(\ung)}
\def\ox{\overline{x}}
\def\uqgh{{{U_{q}(\ungh)}}}
\def\uepfinh{{U_\epsilon^{\roman{fin}}(\hat\ung)}}

\def\calp{{{{\Cal{P}}}}}
\def\calpd{\Bbb D}
\def\ot{{{\otimes}}}

\def\calp{\Cal P}

\magnification 1200
\input amstex
\documentstyle{amsppt}
\NoBlackBoxes
\nologo
\document
\centerline{\bf{{Quantum Affine Algebras at Roots of Unity }}} 
\vskip 36pt
\centerline{Vyjayanthi Chari{\footnote{Partially supported by NATO  and the ESPRC (GR/K65812)}} and Andrew Pressley{\footnote{Partially supported by NATO and the EPSRC (GR/L26216)}}}

\vskip80pt\centerline{ABSTRACT}

\noindent
{\eightpoint{
Let $\uqgh$ be the quantized universal enveloping algebra of the affine Lie algebra $\hat\ung$ associated to a finite-dimensional complex simple Lie algebra $\ung$, and let $\uqresh$ be the $\Bbb C[q,q^{-1}]$-subalgebra of $\uqgh$ generated by the $q$-divided powers of the Chevalley generators. Let $\uepresh$ be the Hopf algebra obtained from $\uqresh$ by specialising $q$ to a non-zero complex number $\epsilon$ of odd order. We classify the finite-dimensional irreducible representations of $\uepresh$ in terms of highest weights. We also give a \lq factorisation\rq\ theorem for such representations: namely, any finite-dimensional irreducible representation of $\uepresh$ is isomorphic to a tensor product of two representations, one factor being the pull-back of a representation of $\hat\ung$ by Lusztig's Frobenius homomorphism $\hat{\roman{Fr}}_\epsilon:\uepresh\to U(\hat\ung)$, the other factor being an irreducible representation of the Frobenius kernel. Finally, we give a conc!
!
!
rete construction of all of the finite-dimensional irreducible representations of $U_\epsilon^{\roman{res}}(\hat{sl}_2)$. The proofs make use of several interesting new identities in $\uqgh$.
}}

\topspace{1cm}

\pagebreak

\noindent{\bf 0. Introduction}
\vskip12pt\noindent
Let $U_q(\ung)$ be the Drinfel'd--Jimbo quantum group associated to a symmetrizable Kac--Moody algebra $\ung$. Thus, $U_q(\ung)$ is a Hopf algebra over the field $\Bbb C(q)$ of rational functions of an indeterminate $q$, and is defined by certain generators and relations (which are written down in Proposition 1.1 below for the cases of interest in this paper). Roughly speaking, one thinks of $U_q(\ung)$ as a family of Hopf algebras depending on a \lq parameter\rq\ $q$. To make this precise, one constructs a \lq specialisation\rq\ $U_\epsilon(\ung)$ of $U_q(\ung)$, for a non-zero complex number $\epsilon$, as follows. If $\epsilon\in\Bbb C$ is transcendental, define $U_\epsilon(\ung)=U_q(\ung)\otimes_{\Bbb C(q)}\Bbb C$, via the algebra homomorphism $\Bbb C(q)\to\Bbb C$ that takes $q$ to $\epsilon$. If, on the other hand, $\epsilon$ is algebraic, the latter homomorphism does not exist, and one proceeds by first constructing a $\Bbb C[q,q^{-1}]$-form of $U_q(\ung)$, i.e. a $\Bbb !
!
!
C[q,q^{-1}]$-(Hopf) subalgebra $\tilde{U}_q(\ung)$ of $U_q(\ung)$ such that $U_q(\ung)=\tilde{U}_q(\ung)\otimes_{\Bbb C[q,q^{-1}]}\Bbb C(q)$. Then one defines $U_\epsilon(\ung)=\tilde{U}_q(\ung)\otimes_{\Bbb C[q,q^{-1}]}\Bbb C$, via the algebra homomorphism $\Bbb C[q,q^{-1}]\to\Bbb C$ that takes $q$ to $\epsilon$. 

Two such $\Bbb C[q,q^{-1}]$-forms have been studied in the literature. They lead to the same algebra $U_\epsilon(\ung)$ when $\epsilon$ is not a root of unity, but different algebras, with very different representation theories, when $\epsilon$ is a root of unity. 

In the \lq non-restricted\rq\ form, one takes $\tilde{U}_q(\ung)$ to be the $\Bbb C[q,q^{-1}]$-subalgebra of ${U}_q(\ung)$ generated by the Chevalley generators $e_i$, $f_i$ of $U_q(\ung)$. The finite-dimensional representations of the non-restricted specialisation $U_\epsilon(\ung)$ have been studied by De Concini, Kac and Procesi [10], [11], when $\ung$ is finite-dimensional, and by Beck and Kac [2], when $\ung$ is (untwisted) affine. 

In the \lq restricted\rq\  form, one takes $\tilde{U}_q(\ung)$ to be the $\Bbb C[q,q^{-1}]$-subalgebra of ${U}_q(\ung)$ generated by the divided powers $e_i^r/[r]_q!$ and $f_i^r/[r]_q!$, for all $r>0$, where $[r]_q!$ denotes a $q$-factorial. When $\ung$ is finite-dimensional and $\epsilon$ is a root of unity, the structure and representation theory of the restricted specialisation $U_\epsilon^{\roman{res}}(\ung)$ was worked out by Lusztig (see [5], [15] and the references there). It is the purpose of the present work to study the finite-dimensional representations of $U_\epsilon^{\roman{res}}(\ung)$ when $\epsilon$ is a root of unity and $\ung$ is (untwisted) affine.

Part of our work may be regarded as the quantum analogue of Garland's paper [13]. A crucial role is played in [13] by a certain identity (Lemma 7.5) that is needed to prove a suitable triangular decomposition of the restricted form of $U(\hat\ung)$ (analogous to $\uqresh$). The proof of the analogue of Garland's lemma in the quantum case (Lemma 5.1 below) is, however, more difficult than in the classical situation because the generators of the \lq positive part\rq\ of $U_\epsilon^{\roman{res}}(\hat\ung)$ do not commute, whereas their classical analogues do commute. Moreover, Garland makes use of a natural derivation of $U(\hat\ung)$ which turns out to have no straightforward analogue in the quantum situation. Lemma 5.1 is one of several interesting new identities in $\uqresh$ that we obtain below. One of these (equation (19) below) shows an unexpected (and as yet unexplained) connection between $U_q^{{\roman{res}}}(\hat{sl}_2)$ and Young diagrams. (This relation is invisible i!
!
!
n the classical situation considered by Garland.)

Once the triangular decomposition of $\uepresh$ is available, we are able to give, in Theorem 8.2, an abstract highest weight description of its finite-dimensional irreducible representations. It turns out that there is a natural one-to-one correspondence between the finite-dimensional irreducible representations of $U_\epsilon^{\roman{res}}(\hat{\ung})$ when $\epsilon$ is a root of unity, and those of $U_\epsilon(\hat\ung)$ when $\epsilon$ is not a root of unity, although corresponding representations have different dimensions, in general (the representations of $U_\epsilon(\hat\ung)$ when $\epsilon$ is not a root of unity were classified in [4], [5] and [6]). This is exactly parallel to the situation for $U_\epsilon^{\roman{res}}(\ung)$ when $\dim(\ung)<\infty$, where the finite-dimensional irreducible representations are parametrised by dominant integral weights whether or not $\epsilon$ is a root of unity (see [5], Chapter 11, and [15]).

One of the most important results about the finite-dimensional irreducible representations $V$ of $U_\epsilon^{\roman{res}}(\ung)$ when $\dim(\ung)<\infty$ and $\epsilon$ is a root of unity of order $\ell$ asserts that $V$ factorises into the tensor product of a representation whose highest weight is divisible by $\ell$ and one whose highest weight is \lq less than $\ell$\rq, in the sense that its value on every simple coroot of $\ung$ is less than $\ell$ (see [5] and [15]). In Section 9, we prove an analogue of this result for $U_\epsilon^{\roman{res}}(\hat{\ung})$ (Theorem 9.1).

It is well known that there is a close relationship between finite-dimensional representations of quantum affine algebras and affine Toda theories (see [7] and the references there). The value of $\epsilon$ is determined by the \lq coupling constant\rq\ of the associated theory. Since the representation theory of $\uepresh$ depends crucially on whether $\epsilon$ is a root of unity or not, one would expect that the structure of affine Toda theories will be different at certain special values of the coupling constant. This does indeed appear to be the case (T. J. Hollowood, private communication), but we shall leave further discussion of this point to another place.
\vskip24pt
\noindent{\bf 1. Preliminaries}
\vskip12pt\noindent
In this section and the next, we recall certain facts about $\ung$ and $\ungh$ and their associated quantum groups that will be needed later. See [1], [5] and [15] for further details.

Let $(a_{ij})_{i,j\in I}$ be the Cartan matrix of the finite-dimensional complex simple Lie algebra $\ung$, let $\hat I=I\coprod\{0\}$, and let $(a_{ij})_{i,j\in\hat I}$ be the generalised Cartan matrix of the untwisted affine Lie algebra $\ungh$ of $\ung$. Let $(d_i)_{i\in\hat I}$ be the coprime positive integers such that the matrix $(d_ia_{ij})_{i,j\in\hat I}$ is symmetric. 

Let $\check P$ be the lattice over $\Bbb Z$ with basis $(\check{\lambda}_i)_{i\in I}$, let $\check{\alpha}_j=\sum_{i\in I}a_{ji}\check{\lambda}_i$ ($j\in I$), and let $\check Q=\sum_{i\in I}\Bbb Z \check{\alpha}_i\subseteq\check P$. The root lattice $Q={\roman{Hom}}_{\Bbb Z}(\check{P},\Bbb Z)$ has basis the simple roots $(\alpha_i)_{i\in I}$ of $\ung$, where $\alpha_i(\check{\lambda}_j)=\delta_{ij}$. Similarly, the weight lattice $P={\roman{Hom}}_{\Bbb Z}(\check{Q},\Bbb Z)$ has basis the fundamental weights $(\lambda_i)_{i\in I}$ of $\ung$, where $\lambda_i(\check{\alpha}_j)=\delta_{ij}$. Let $P^+=$

\noindent$\{\lambda\in P\,\mid\, \lambda(\check{\alpha}_i)\ge 0\ \text{for all $i\in I$}\}$, and let $Q^+=\sum_{i\in I}\Bbb N.\alpha_i$.

For $i\in I$, define $s_i:\check P\to\check P$ by $s_i(x)=x-\alpha_i(x)\check{\alpha}_i$, and let $W$ be the group of automorphisms of $\check P$ generated by the $s_i$. Then, $W$ acts on $Q$ by $s_i(\xi)=\xi-\xi(\check{\alpha}_i)\alpha_i$, for $\xi\in Q$. The root and coroot systems of $\ung$ are given, respectively, by $R=\bigcup_{i\in I}W\alpha_i$, $\check{R}=\bigcup_{i\in I}W\check{\alpha}_i$. There is a partial order on $R$ such that $\alpha\le \beta$ if and only if $\beta-\alpha$ is a linear combination of the $\alpha_i$ (${i\in I})$ with non-negative integer coefficients. The correspondence $\alpha_i\leftrightarrow\check{\alpha}_i$ extends to a $W$-invariant correspondence $R\leftrightarrow \check{R}$, written $\alpha\leftrightarrow\check{\alpha}$, such that $\alpha(\check{\alpha})=2$ for all $\alpha\in R$. Define $s_\alpha:\check P\to\check P$ by $s_\alpha(x)=x-\alpha(x)\check{\alpha}$, for all $\alpha\in R$. 

Let $\hat W=W\tilde{\times}\check P$ be the semi-direct product group defined using the $W$ action on $\check P$. For $s\in W$, write $s$ for $(s,0)\in\hat W$, and for $x\in\check P$ write $x$ for $(1,x)\in\hat W$, where $1$ is the identity element of $W$. Let $\theta$ be the highest root of $\ung$ with respect to $\le$, and write $s_0$ for $(s_\theta,\check\theta)\in\hat W$. Let $\tilde W$ be the (normal) subgroup of $\hat W$ generated by the $s_i$ for $i\in\hat I$, and let $\Tau=\hat{W}/\tilde{W}$. Then, $\Tau$ is a finite group isomorphic to a subgroup of the group of diagram automorphisms of $\ungh$, i.e. the bijections $\tau:\hat I\to\hat I$ such that $a_{\tau(i)\tau(j)}=a_{ij}$ for all $i,j\in\hat I$. Moreover, there is an isomorphism of groups $\hat W\cong\Tau\tilde{\times}\tilde W$, where the semi-direct product is defined using the action of $\Tau$ in $\tilde{W}$ given by $\tau.s_i=s_{\tau(i)}$ (see [1]). If $w\in\hat W$, a reduced expression for $w$ is an expression !
!
!
$w=\tau s_{i_1}s_{i_2}\ldots s_{i_n}$ with $\tau\in\Tau$, $i_1,i_2,\ldots,i_n\in\hat I$ and $n$ minimal. 

The universal enveloping algebra $\overline{U}(\ung)$ (resp. $\overline{U}(\ungh)$) is the associative algebra over $\Bbb C$ with generators $\overline{e}_i^\pm$, $\overline{h}_i$ for $i\in I$ (resp. $i\in\hat I$) and defining relations
$$\align
\overline{h}_i\overline{h}_j=\overline{h}_j\overline{h}_i,\ \ \ &\oh_i\oe_j^\pm-\oe_j^\pm\overline{h}_i=\pm a_{ij}\oe_j^\pm,\\
\oe_i^+\oe_j^--\oe_j^-\oe_i^+&=\delta_{ij}\oh_i,\\
\sum_{r=0}^{1-a_{ij}}(-1)^r\left({1-a_{ij}}\atop r\right)&(\oe_i^\pm)^r\oe_j^\pm(\oe_i^\pm)^{1-a_{ij}-r}=0\ \ \ \ \ \text{if $i\ne j$},\endalign$$
where $i,j\in I$ (resp. $i,j\in\hat I$). 

Let $q$ be an indeterminate, let $\Bbb C(q)$ be the field of rational functions of $q$ with complex coefficients, and let $\Bbb C[q,q^{-1}]$ be the ring of Laurent polynomials. For $r,n\in\Bbb N$, $n\ge r$, define
$$\align [n]_q & =\frac{q^n -q^{-n}}{q -q^{-1}},\\
[n]_q! &=[n]_q[n-1]_q\ldots [2]_q[1]_q,\\
\left[{n\atop r}\right]_q &= \frac{[n]_q!}{[r]_q![n-r]_q!}.\endalign$$
Then $\left[{n\atop r}\right]_q\in\Bbb C[q,q^{-1}]$ for all $n\ge r\ge 0$. One has the identity 
$$\sum_{r=0}^n(-1)^rq^{r(n-1)}\left[{n\atop r}\right]_q=0\tag1$$
for all $n\ge 0$. Let $q_i=q^{d_i}$ for $i\in\hat I$.

\proclaim{Proposition 1.1} There is a Hopf algebra $\uqgh$ over $\Bbb C(q)$ which is generated as an algebra by elements $e_i^\pm$, $k_i^{{}\pm 1}$ ($i\in\hat I$), with the following defining relations:
$$\align 
k_ik_i^{-1}=k_i^{-1}k_i&=1,\ \ \ \ k_ik_j=k_jk_i,\\
k_ie_j^\pm k_i^{-1}&=q_i^{\pm a_{ij}}e_j^\pm,\\
[e_i^+,e_j^-]&=\delta_{ij}\frac{k_i-k_i^{-1}}{q_i-q_i^{-1}},\\
\sum_{r=0}^{1-a_{ij}}(-1)^r\left[{1-a_{ij}}\atop r\right]_{q_i}&(e_i^\pm)^re_j^\pm(e_i^\pm)^{1-a_{ij}-r}=0\ \ \ \ \ \text{if $i\ne j$}.
\endalign$$
The comultiplication of $\uqgh$ is given on generators by
$$\Delta(e_i^+)=e_i^+\ot k_i+1\ot e_i^+,\ \ 
\Delta(e_i^-)=e_i^-\ot 1+k_i^{-1}\ot e_i^-,\ \ 
\Delta(k_i)=k_i\ot k_i,$$
for $i\in\hat I$.

Restricting the indices $i,j$ to $I$, one obtains a Hopf algebra $\uqg$. 
\endproclaim

Define $\uqg^\pm$ (resp. $\uqg^0$) to be the $\Bbb C(q)$-subalgebra of $\uqg$ generated by the $e_i^\pm$ (resp. by the $k_i^{\pm 1}$) for all $i\in I$. 

There is a natural injective homomorphism of $\Bbb C(q)$-Hopf algebras $\uqg\to\uqgh$ that takes $e_i^\pm$ to $e_i^\pm$ and $k_i$ to $k_i$, for all $i\in I$. In particular, any representation of $\uqgh$ can be regarded as a representation of $\uqg$. 

It is convenient to use the following notation:
$$(e_i^\pm)^{(r)}=\frac{(e_i^\pm)^r}{[r]_{q_i}!}.$$

Let $T_i$ ($i\in\hat I$) be the $\Bbb C(q)$-algebra automorphisms of $\uqgh$ defined by Lusztig [15] (our notational conventions differ slightly from his):
$$\align
T_i(e_i^+)&=-e_i^-k_i,\ \ T_i(e_i^-)=-k_i^{-1}e_i^+,\ \ T_i(k_j)=k_i^{-a_{ij}}k_j,\\
T_i(e_j^+)&=\sum_{r=0}^{-a_{ij}}(-1)^{r-a_{ij}}q_i^{-r}(e_i^+)^{(-a_{ij}-r)}e_j^+(e_i^+)^{(r)}\ \ \text{if $i\ne j$},\\
T_i(e_j^-)&=\sum_{r=0}^{-a_{ij}}(-1)^{r-a_{ij}}q_i^{r}(e_i^-)^{(r)}e_j^-(e_i^-)^{(-a_{ij}-r)}\ \ \text{if $i\ne j$}.\endalign$$
If $i\in I$, $T_i$ induces an automorphism of the subalgebra $\uqg$ of $\uqgh$, also denoted by $T_i$. 

There is, of course, a classical analogue $\overline{T}_i$ of $T_i$ for all $i\in \hat I$, given by the above formulas with $e_i$ replaced by $\overline{e}_i$ and $k_i$ replaced by $1$, and with $\overline{T}_i(\oh_j)=\oh_j-a_{ij}\oh_i$. The $\overline{T}_i$ are Hopf algebra automorphisms of $\overline{U}(\ungh)$. 

The finite group $\Cal T$ acts as $\Bbb C(q)$-Hopf algebra automorphisms of $\uqgh$ by
$$\tau.e_i^\pm=e_{\tau(i)}^\pm,\ \ \ \tau.k_i=k_{\tau(i)}, \ \ \ \text{for all $i\in\hat I$}.$$
Similarly, $\Tau$ acts as Hopf algebra automorphisms of $\overline{U}(\ungh)$ by
$$\tau.\oe_i^\pm=\oe_{\tau(i)}^\pm,\ \ \ \tau.\oh_i=\oh_{\tau(i)}, \ \ \ \text{for all $i\in\hat I$}.$$

If $w\in W$ has a reduced expression $w=\tau s_{i_1}\ldots s_{i_n}$, let $T_w$ be the $\Bbb C(q)$-algebra automorphism of $\uqgh$ given by $T_w=\tau T_{i_1}\ldots T_{i_n}$. Then, $T_w$ depends only on $w$, and not on the choice of its reduced expression [1]. In particular, for any $i\in I$, we have a well defined $\Bbb C(q)$-algebra automorphism $T_{\check{\lambda}_i}$ of $\uqgh$. Similarly, one has Hopf algebra automorphisms $\overline{T}_{\check{\lambda}_i}$ of $\overline{U}(\ungh)$. 

Let $\uqres$ (resp. $\uqresh$) be the $\Bbb C[q,q^{-1}]$-subalgebra of $\uqg$ 
(resp. $\uqgh$) generated by the $k_i^{\pm 1}$ and the $(e_i^\pm)^{(r)}$ for all $i\in I$, $r\in\Bbb N$ (resp. for all $i\in\hat I$, $r\in\Bbb N$). Then, $\uqres$ (resp. $\uqresh$) is a $\Bbb C[q,q^{-1}]$-Hopf subalgebra of $\uqg$ (resp. $\uqgh$). Moreover, $\uqres$ (resp. $\uqresh$) is preserved by the automorphism $T_i$ for all $i\in I$ (resp. for all $i\in \hat I$). In fact [15], if $i,j\in\hat I$, $n\in\Bbb N$,
$$\align
T_i((e_i^+)^{(n)})&=(-1)^nq_i^{-n(n-1)}(e_i^-)^{(n)}k_i^n,\ \ T_i((e_i^-)^{(n)})=(-1)^nq^{n(n-1)}k_i^{-n}(e_i^+)^{(n)},\\
T_i((e_j^+)^{(n)})&=\sum_{r=0}^{-na_{ij}}(-1)^{r-na_{ij}}q_i^{-r}
(e_i^+)^{(-na_{ij}-r)}(e_j^+)^{(n)}(e_i^+)^{(r)}\ \ \text{if $i\ne j$},\\
T_i((e_j^-)^{(n)})&=\sum_{r=0}^{-na_{ij}}(-1)^{r-na_{ij}}q_i^{r}(e_i^-)^{(r)}
(e_j^-)^{(n)}(e_i^-)^{(-na_{ij}-r)}\ \ \text{if $i\ne j$}.\endalign$$
The action of $\Tau$ obviously preserves $\uqresh$. 

For $r\in\Bbb N$, $n\in\Bbb Z$, $i\in I$, define
$$\left[{k_i;n}\atop r\right]=\prod_{s=1}^r\frac{k_iq_i^{n-s+1}-k_i^{-1}q_i^{-n+s-1}}
{q_i^s-q_i^{-s}}.$$
The importance of these elements stems from the identity
$$(e_i^+)^{(r)}(e_i^-)^{(s)}=\sum_{t=0}^{{\text{min}}(r,s)}(e_i^-)^{(s-t)}
\left[{k_i;2t-r-s}\atop t\right](e_i^+)^{(r-t)},\tag2$$
for all $r,s\in\Bbb N$, $i\in I$. It follows from this identity that 
$\left[{k_i;n}\atop r\right]\in\uqres$ for all $r\in\Bbb N$, $n\in\Bbb Z$, $i\in I$.

Let $U_q^{{\roman{res}}}(\ung)^\pm$ be the subalgebra of $\uqres$ generated by the $(e_i^\pm)^{(r)}$ for all $i\in I$, $r\in\Bbb N$, and let $U_q^{{\roman{res}}}(\ung)^0$ be the subalgebra generated by the $k_i^{\pm 1}$ and the $\left[{k_i;n}\atop r\right]$, for all $r\in\Bbb N$, $n\in\Bbb Z$, $i\in I$. Multiplication defines an isomorphism of $\Bbb C(q)$-vector spaces
$$U_q^{{\roman{res}}}(\ung)^-\otimes U_q^{{\roman{res}}}(\ung)^0\otimes U_q^{{\roman{res}}}(\ung)^+\to\uqres.$$
The appropriate definitions of $U_q^{{\roman{res}}}(\hat\ung)^\pm$ and $U_q^{{\roman{res}}}(\hat\ung)^0$ will be given at the end of Section 3.

Let $\ell$ be an odd integer with $\ell\ge 3$, and assume that $\ell$ is not divisible by $3$ if $\ung$ is of type $G_2$. Let $\epsilon\in\Bbb C$ be a primitive $\ell$th root of unity, and set
$$\uepres=\uqres\otimes_{\Bbb C[q,q^{-1}]}\Bbb C,\ \ \ 
\uepresh=\uqresh\otimes_{\Bbb C[q,q^{-1}]}\Bbb C,$$
via the algebra homomorphism $\Bbb C[q,q^{-1}]\to\Bbb C$ that takes $q$ to $\epsilon$. If $x$ is any element of $\uqres$ (resp. $\uqresh$), we denote the corresponding element of $\uepres$ (resp. $\uepresh$) also by $x$. 

Define 
$U_\epsilon^{{\roman{res}}}(\ung)^\pm$ and $U_\epsilon^{{\roman{res}}}(\ung)^0$
in the obvious way. It is known that 
$U_\epsilon^{{\roman{res}}}(\ung)^0$ is generated by the $k_i^{\pm 1}$ and the $\left[{k_i;0}\atop{\ell}\right]$, for $i\in I$. Multiplication defines an isomorphism of $\Bbb C$-vector spaces
$$U_\epsilon^{{\roman{res}}}(\ung)^-\otimes U_\epsilon^{{\roman{res}}}(\ung)^0\otimes U_\epsilon^{{\roman{res}}}(\ung)^+\to\uepres.$$
Let $U_\epsilon^{{\roman{fin}}}(\ung)$ be the subalgebra of $U_\epsilon^{{\roman{res}}}(\ung)$ generated by the $e_i^\pm$ for $i\in I$ and 
$U_\epsilon^{{\roman{res}}}(\ung)^0$.

We shall make use of the characteristic zero Frobenius homomorphisms defined by Lusztig (see [15], Chapter 35). Namely, there exist homomorphisms of Hopf algebras
$${\roman{Fr}}_\epsilon:U_\epsilon^{{\roman{res}}}(\ung)\to\overline{U}(\ung),\ \ \ \ \hat{\roman{Fr}}_\epsilon:\uepresh\to\overline{U}(\hat{\ung})$$
such that, for $i\in I$,
$${\roman{Fr}}_\epsilon(k_i)=1,\ \ {\roman{Fr}}_\epsilon((e_i^\pm)^{(r)})=\cases \frac{(\overline{e}_i^\pm)^{r/\ell}}{(r/\ell)!} & \text{if $\ell$ divides $r$},\\
0 & \text{otherwise.}
\endcases$$
and for $i\in\hat I$,
$$\hat{\roman{Fr}}_\epsilon(k_i)=1,\ \ \hat{\roman{Fr}}_\epsilon(({e}_i^\pm)^{(r)})=\cases \frac{(\overline{e}_i^\pm)^{r/\ell}}{(r/\ell)!} & \text{if $\ell$ divides $r$},\\
0 & \text{otherwise.}
\endcases$$

To discuss finite-dimensional representations of $\uqresh$, or of its specialization $\uepresh$, it is advantageous to use another realization of $\uqgh$, due to Drinfel'd [12], Beck [1] and Jing [14].
\proclaim{Theorem 1.2} There is an isomorphism of $\Bbb C(q)$-Hopf algebras from $\uqgh$ to the algebra with generators $x_{i,r}^{{}\pm{}}$ ($i\in I$, $r\in\Bbb Z$), $k_i^{{}\pm 1}$ ($i\in I$), $h_{i,r}$ ($i\in I$, $r\in \Bbb Z\backslash\{0\}$) and $c^{{}\pm{1/2}}$, and the following defining relations:
$$\align
c^{{}\pm{1/2}}\ &\text{are central,}\\
k_ik_i^{-1} = k_i^{-1}k_i =1,\;\; &c^{1/2}c^{-1/2} =c^{-1/2}c^{1/2} =1,\\
k_ik_j =k_jk_i,\;\; &k_ih_{j,r} =h_{j,r}k_i,\\
k_ix_{j,r}k_i^{-1} &= q_i^{{}\pm a_{ij}}x_{j,r}^{{}\pm{}},\\
[h_{i,r},h_{j,s}]&=\delta_{r,-s}\frac1{r}[ra_{ij}]_{q_i}\frac{c^r-c^{-r}}
{q_j-q_j^{-1}},\\
[h_{i,r} , x_{j,s}^{{}\pm{}}] &= \pm\frac1r[ra_{ij}]_{q_i}c^{{}\mp {|r|/2}}x_{j,r+s}^{{}\pm{}},\\
x_{i,r+1}^{{}\pm{}}x_{j,s}^{{}\pm{}} -q_i^{{}\pm a_{ij}}x_{j,s}^{{}\pm{}}x_{i,r+1}^{{}\pm{}} &=q_i^{{}\pm a_{ij}}x_{i,r}^{{}\pm{}}x_{j,s+1}^{{}\pm{}} -x_{j,s+1}^{{}\pm{}}x_{i,r}^{{}\pm{}},\tag3\\
[x_{i,r}^+ , x_{j,s}^-]=\delta_{ij} & \frac{ c^{(r-s)/2}\psi_{i,r+s}^+ - c^{-(r-s)/2} \psi_{i,r+s}^-}{q_i - q_i^{-1}},\\
\sum_{\pi\in\Sigma_m}\sum_{k=0}^m(-1)^k\left[{m\atop k}\right]_{q_i} x_{i, r_{\pi(1)}}^{{}\pm{}}\ldots x_{i,r_{\pi(k)}}^{{}\pm{}} & x_{j,s}^{{}\pm{}}
 x_{i, r_{\pi(k+1)}}^{{}\pm{}}\ldots x_{i,r_{\pi(m)}}^{{}\pm{}} =0,\ \ \text{if $i\ne j$},
\endalign$$
for all sequences of integers $r_1,\ldots, r_m$, where $m =1-a_{ij}$, $\Sigma_m$ is the symmetric group on $m$ letters, and the $\psi_{i,r}^{{}\pm{}}$ are determined by equating powers of $u$ in the formal power series 
$$\sum_{r=0}^{\infty}\psi_{i,\pm r}^{{}\pm{}}u^{{}\pm r} = k_i^{{}\pm 1} exp\left(\pm(q_i-q_i^{-1})\sum_{s=1}^{\infty}h_{i,\pm s} u^{{}\pm s}\right).$$

The isomorphism is given by 
$$x_{i,r}^\pm=o(i)^rT_{\check{\lambda}_i}^{\mp r}(e_i^\pm),$$
where $o:I\to\{\pm 1\}$ is a map such that $o(i)=-o(j)$ whenever $a_{ij}<0$ (it is clear that there are exactly two possible choices for $o$).
\endproclaim

The classical analogues of the generators $x_{i,r}^\pm$ and $h_{i,r}$ are given by
$$\ox_{i,r}^\pm=o(i)^r\overline{T}_{\check{\lambda}_i}^{\mp r}(\oe_i^\pm),
\ \ \ \ \oh_{i,r}=[\ox_{i,r}^+,\ox_{i,0}^-].$$
The latter elements, together with a central element $\overline{c}$, generate $\overline{U}(\hat\ung)$ and the defining relations are obtained from those in 1.2 by interpreting $c$ as $q^{\overline c}$, $k_i$ as $q_i^{\oh_{i,0}}$, $\psi_{i,r}^+$ (resp. $\psi_{i,r}^-$) as $(q_i-q_i^{-1})\oh_{i,r}$ if $r>0$ (resp. if $r<0$), and letting $q$ tend to one.

For each $i\in I$, there is an injective homomorphism of $\Bbb C(q)$-algebras $\varphi_i:U_{q_i}(\hat{sl}_2)\to \uqgh$ that takes $x_r^\pm$ to $x_{i,r}^\pm$ and $\psi_r^\pm$ to $\psi_{i,r}^\pm$, for all $r\in\Bbb Z$ (as usual, we drop the (unique) Dynkin diagram subscript when dealing with $U_q(\hat{sl}_2)$) (see [1], Proposition 3.8). It follows immediately from 1.2 that $U_{q_i}^{\roman{res}}(\hat{sl}_2)$ is generated as a $\Bbb C[q,q^{-1}]$-algebra by $c^{\pm 1/2}$, $k^{\pm 1}$ and the $(x_n^\pm)^{(r)}$ for $n\in\Bbb Z$, $r\in\Bbb N$ (for $e_0^+=x_1^-$ and $e_0^1=c^{-1}kx_{-1}^+$). Since the $T_{\check{\lambda}_i}$ preserve $\uqresh$, it is also clear that $(x_{i,n}^\pm)^{(r)}\in\uqresh$ for all $i\in I$, $n\in\Bbb Z$, $r\in\Bbb N$. It follows that $\varphi_i(U_{q_i}^{\roman{res}}(\hat{sl}_2))\subseteq\uqresh$, and hence that $\varphi_i$ induces a $\Bbb C$-algebra homomorphism $\varphi_i:U_{\epsilon_i}^{\roman{res}}(\hat{sl}_2)\to\uepresh$.
We have that $T_{\check{\lambda}_j}$ fixes $\varphi_i(U_{q_i}(\hat{sl}_2))$ pointwise if $i,j\in I$, $i\ne j$ ([1], Corollary 3.2). We denote the obvious classical analogue of $\varphi_i$ by $\overline{\varphi}_i:\overline{U}(\hat{sl}_2)\to\overline{U}(\ungh)$.

Finally, we shall need the following automorphisms of $\uqgh$:

\proclaim{Proposition 1.3}  (a) There exists a $\Bbb C(q)$-algebra anti-automorphism $\Phi$ of $\uqgh$ such that, for all $i\in I, r\in\Bbb Z$,
$$\Phi(x_{i,r}^\pm)=x_{i,r}^\mp,\ \ \Phi(c^{1/2})=c^{-1/2},\ \  \Phi(\psi_{i,r}^\pm)=\psi_{i,r}^\pm,\ \ \Phi(h_{i,r})=h_{i,r}.$$

\noindent (b) There exists a $\Bbb C$-algebra and coalgebra anti-automorphism $\Omega$ of $\uqgh$ such that $\Omega(q)=q^{-1}$ and, for all $i\in I, r\in\Bbb Z$,
$$\Omega(x_{i,r}^\pm)=x_{i,-r}^\mp,\ \ \Omega(\psi_{i,r}^\pm)=\psi_{i,-r}^\mp,\ \ 
\Omega(h_{i,r}) =h_{i,-r},\ \  \Omega(c^{1/2})=c^{-1/2}.$$
Moreover, $\Omega$ and $\Phi$ each preserve $\uqresh$, commute with $\Tau$, and we have, for all $i\in I$,
$$\Omega T_i=T_i\Omega,\ \ \ \Phi T_i=T_i^{-1}\Phi,\ \ \ \Phi\Omega=\Omega\Phi.$$
\endproclaim
\demo{Proof} The existence of $\Phi$ and $\Omega$ can be verified directly by using 1.2. It is obvious that $\Phi$ and $\Omega$ commute with $\Tau$. The relations between $\Omega$, $\Phi$ and the $T_i$ are checked by verifying agreement on the generators $e_j^\pm$,  $k_j^{\pm 1}$ ($j\in\hat I$) of $\uqgh$, on which they are given by
$$\align
\Phi(e_j^\pm)=e_j^\mp,\ \ \ \ &\Phi(k_j)=k_j,\\
\Omega(e_j^\pm)=e_j^\mp,\ \ \ \ &\Omega(k_j)=k_j^{-1}.
\endalign$$
These formulas obviously imply that $\Phi$ and $\Omega$ preserve $\uqresh$. The fact that $\Omega$ is a coalgebra anti-automorphism can also be checked on the generators $e_j^\pm,k_j^{\pm 1}$ ($j\in\hat I$). \qed\enddemo

\vskip24pt\noindent
{\bf 2. Representation Theory of $U_q$}
\vskip 12pt\noindent
A representation $V$ of $\uqg$ is said to be of type I if, for all $i\in I$, $k_i$ acts semisimply on $V$ with eigenvalues in $q_i^{\Bbb Z}$. Any finite-dimensional irreducible representation of $\uqg$ can be obtained from a type I representation by tensoring with a one-dimensional representation that takes each $e_i^\pm$ to zero and each $k_i$ to $\pm 1$.

If $\lambda\in P$, the weight space $V_\lambda$ of a representation $V$ of $\uqg$ is defined by
$$V_\lambda=\left\{v\in V \bigm| k_i.v=q_i^{n_i}v\right\},$$
where $n_i=\lambda(\check{\alpha}_i)$.
We have $$e_i^{{}\pm{}}.V_\lambda\subseteq V_{\lambda\pm\alpha_i}.$$

A vector $v$ in a type I representation $V$ of $\uqg$ is said to be a highest weight vector if there exists $\lambda\in P$ such that 
$$v\in V_\lambda\ \ \text{and}\ \ e_i^+.v=0\ \ \ {\text{for all}} \  r\in\Bbb N, i\in I.$$
If, in addition, $V=\uqg.v$, then $V$ is said to be a highest weight representation with highest weight $\lambda$.

For any $\lambda\in P$, there exists, up to isomorphism, a unique irreducible representation $V(\lambda)$ of $\uqg$ with highest weight $\lambda$. We have
$$V(\lambda)=\bigoplus_{\mu\le\lambda}V(\lambda)_\mu.$$
Every finite-dimensional irreducible type I representation of $\uqg$ is isomorphic to $V(\lambda)$ for some (unique) $\lambda\in P^+$. The character of $V(\lambda)$, and in particular its dimension, is the same as that of the irreducible representation $\overline{V}(\lambda)$ of $\ung$ with the same highest weight $\lambda$. 

Turning now to $\uqgh$, we say that a representation $V$ of $\uqgh$ is type I if it is type I as a representation of $\uqg$ and $c^{1/2}$ acts as one on $V$. A vector $v$ in a type I representation $V$ of $\uqgh$ is said to be a highest weight vector if $v$ is annihilated by $x_{i,r}^+$, and is an eigenvector of $\psi_{i,r}^\pm$, for all $i\in I,r\in\Bbb Z$. If, in addition, $V=\uqgh.v$, then $V$ is said to be a highest weight representation of $\uqgh$. One shows by the usual Verma module arguments that, for any set of scalars $\Psi_{i,r}^\pm$ such that $\Psi_{i,0}^+\Psi_{i,0}^-=1$, there is an irreducible highest weight representation $V$ of $\uqgh$ with highest weight vector $v$ such that $\psi_{i,r}^\pm.v=\Psi_{i,r}^\pm v$ for all $i\in I,r\in\Bbb Z$ (compare the proof of 7.3 below). Moreover, $V$ is uniquely determined by the scalars $\Psi_{i,r}^\pm$, up to isomorphism.

\proclaim{Theorem 2.1} The irreducible highest weight representation of $\uqgh$ determined by the scalars $\Psi_{i,r}^\pm$ is finite-dimensional if and only if there exist polynomials $P_i\in\Bbb C(q)[u]$, for $i\in I$, such that $P_i(0)\ne 0$ and 
$$\sum_{r=0}^\infty\Psi_{i,r}^+u^r=q_i^{{\roman{deg}}(P_i)}
\frac{P_i(q_i^{-2}u)}{P_i(u)}=\sum_{r=0}^\infty\Psi_{i,-r}^-u^{-r},$$
in the sense that the left- and right-hand sides are the Laurent expansions of the middle term about $u=0$ and $u=\infty$, respectively.\endproclaim

This is proved in [4], [5], [6] and [9]. From now on, we denote by $\Pi_q$ the set of polynomials $P\in\Bbb C(q)[u]$ such that $P(0)\ne 0$, and by $\Pi_q^I$ the set of $I$-tuples of such polynomials. If $\bold P=(P_i)_{i\in I}\in\Pi_q^I$, we denote by $V(\bold P)$ the finite-dimensional irreducible representation of $\uqgh$ determined by $\bold P$ as in 2.1.

If $\bold P=(P_i)_{i\in I}$, $\bold Q=(Q_i)_{i\in I}\in\Pi_q^I$, define $\bold P\otimes\bold Q=(P_iQ_i)_{i\in I}$. The following result is proved in [9] (see 7.4 for the proof of an analogous result):

\proclaim{Proposition 2.2} Let $\bold P,\bold Q\in\Pi_q^I$. If $V(\bold P)\otimes V(\bold Q)$ is irreducible as a representation of $\uqgh$, then it is isomorphic to $V(\bold P\otimes\bold Q)$.\endproclaim

We shall also need the classical analogue of 2.1. Let $\Pi$ the set of polynomials $P\in\Bbb C[u]$ such that $P(0)\ne 0$, and let $\Pi^I$ be the set of $I$-tuples of such polynomials. Then, the classical analogue of 2.1 states that every finite-dimensional irreducible representation $\overline{V}$ of $\overline{U}(\hat\ung)$ is generated by a vector ${v}$ that is annihilated by $\overline c$ and the $\ox_{i,r}^\pm$, and such that $\oh_{i,r}.v={H}_{i,r}v$, where the scalars 
${H}_{i,r}\in\Bbb C$ satisfy
$$\sum_{r=0}^\infty {H}_{i,r}u^r={\roman{deg}}({P}_i)-u\frac{{P}_i'(u)}
{{P}_i(u)},\ \  \sum_{r=0}^\infty {H}_{i,-r}u^{-r}=u\frac{{P}_i'(u)}{{P}_i(u)},$$
for some ${\bold P}=({P}_i)_{i\in I}\in\Pi^I$. We write $\overline{V}$ as $\overline{V}({\bold P})$. The obvious classical analogue of 2.2 holds.

Finite-dimensional irreducible representations of $\uqgh$ can be constructed explicitly when $\ung=sl_2$, by means of the following result. We take $\hat I=\{0,1\}$ when $\ung=sl_2$.

\proclaim{Proposition 2.3} For any non-zero $a\in\Bbb C(q)$, there is an  algebra homomorphism $ev_a:U_q(\hat{sl}_2)\to U_q(sl_2)$ such that
$${\text{ev}}_a(e_0^{{}\pm{}})=q^{{}\pm 1}a^{{}\pm 1}e^{{}\mp{}},\ {\text{ev}}_a(k_0)=k^{-1},\ {\text{ev}}_a(e_1^{{}\pm{}})=e^{{}\pm{}},\ {\text{ev}}_a(k_1) =k.\tag4$$
Moreover, we have, for all $r\in\Bbb Z$,
$${\text{ev}}_a(x_r^+)=q^{-r}a^rk^re^+,\ \ {\text{ev}}_a(x_r^-)=q^{-r}a^re^-k^r.\tag5$$
\endproclaim
See [4], Proposition 5.1, for the proof. If $V$ is a representation of $ U_q(sl_2)$, we denote by $V_a$ the \lq evaluation representation' of $U_q(\hat{sl}_2)$ obtained by pulling back $V$ by ${\roman{ev}}_a$. Let $V(n)$ be the $(n+1)$-dimensional irreducible type I representation of $U_q(sl_2)$.

\proclaim{Proposition 2.4} Let $r,n_1,\ldots,n_r\in\Bbb N$, and let $a_1,\ldots,a_r\in\Bbb C(q)$ be non-zero. Then, the tensor product
$$V(n_1)_{a_1}\otimes\cdots\otimes V(n_r)_{a_r}$$
is irreducible as a representation of $U_q(\hat{sl}_2)$ if and only if, for all $1\le s\ne t\le r$,
$$a_s/a_t\ne q^{\pm(n_s+n_t-2p)}\ \ \ \text{for all $0\le p<{\roman{min}}(n_s,n_t)$}.$$
\endproclaim

This is proved by exactly the same argument as Theorem 5.8 in [4]. 

It is easy to compute that $V(n)_a\cong V(P_{n,a})$, where 
$$P_{n,a}(u)=\prod_{m=1}^n(1-aq^{n-2m+1}u)\tag6$$
(see [4], Corollary 4.2). It follows from 2.2 that, when the conditions in 2.4 are satisfied, 
$$V(n_1)_{a_1}\otimes\cdots\otimes V(n_r)_{a_r}\cong V(P),$$
where $P=\prod_{s=1}^rP_{n_s,a_s}$.
\vskip12pt
The classical analogues of the homomorphisms ${\roman{ev}}_a$ exist for all $\ung$. They are the homomorphisms $\overline{\roman{ev}}_a:{\overline{U}}(\hat\ung)\to\overline{U}(\ung)$, defined for any non-zero $a\in\Bbb C$, such that
$$\overline{\roman{ev}}_a(\overline{e}_0^\pm)=a^{\pm 1}\overline{e}_\theta^\mp,\ \ 
\overline{\roman{ev}}_a(\overline{e}_i^\pm)=\overline{e}_i^\pm,\ \ 
\overline{\roman{ev}}_a(\overline{h}_0)=-\check\theta,\ \ 
\overline{\roman{ev}}_a(\overline{h}_i)=h_i,$$
for $i\in I$. Here, $\overline{e}_\theta^\pm$ are root vectors corresponding to $\pm\theta$, normalised so that $[\overline{e}_\theta^+,\overline{e}_\theta^-]=\check\theta$. Evaluation representations $\overline{V}_a$ of $\overline{U}(\hat\ung)$ are defined in the obvious way. Using this construction, one can describe explicitly the representation $\overline{V}(\bold P)$ for all $\bold P\in\Pi^I$. 

\proclaim{Theorem 2.5} Let $\bold P=(P_i)_{i\in I}\in\Pi^I$, and let $\{a_1^{-1},a_2^{-1},\ldots,a_r^{-1}\}$ be the union of the set of roots of $P_i$ for all $i\in I$. Write
$$P_i(u)=\prod_{s=1}^r(1-a_su)^{n_{i,s}}$$
for some integers $n_{i,s}\ge 0$, and let $\mu_s=\sum_{i\in I}n_{i,s}\lambda_i\in P^+$. Then,
$$\overline{V}(\bold P)\cong \overline{V}(\mu_1)_{a_1}\otimes \overline{V}(\mu_2)_{a_2}\otimes\cdots\otimes \overline{V}(\mu_r)_{a_r}.$$
Conversely, if $\mu_1,\mu_2,\ldots,\mu_r\in P^+$ and $a_1,a_2,\ldots,a_r\in\Bbb C$ are non-zero, the tensor product $\overline{V}(\mu_1)_{a_1}\otimes \overline{V}(\mu_2)_{a_2}\otimes\cdots\otimes \overline{V}(\mu_r)_{a_r}$
is irreducible as a representation of $\hat\ung$ if and only if $a_1,a_2,\ldots,a_r$ are distinct.\endproclaim
\demo{Proof} By using the methods in [3], one shows that every finite-dimensional irreducible representation of $\overline{U}(\hat\ung)$ is isomorphic to a tensor product 
$$\overline{V}(\mu_1)_{a_1}\otimes \overline{V}(\mu_2)_{a_2}\otimes\cdots\otimes \overline{V}(\mu_r)_{a_r}$$
for some $\mu_1,\mu_2,\ldots,\mu_r\in P^+$ and non-zero $a_1,a_2,\ldots,a_r\in\Bbb C$, and that the condition for irreducibility of such a tensor product is as stated in 2.5. By the classical analogue of 2.2, it therefore suffices to prove that the $I$-tuple of polynomials $(P_i)_{i\in I}$ associated to $\overline{V}(\mu_s)_{a_s}$ is given by $P_i(u)=(1-a_su)^{n_{i,s}}$. This is a straightforward computation.\qed\enddemo

\vskip24pt
\noindent{\bf 3. Some Identities}
\vskip12pt\noindent
In this section we establish certain identities in ${U}_q(\hat{sl}_2)$ that will be needed later. By applying the homomorphisms $\varphi_i:U_{q_i}(\hat{sl}_2)\to U_q(\hat{\frak g})$ defined in Section 1, one obtains corresponding identities in $U_q(\hat{\frak g})$. As usual, we drop the subscript $i\in I$ when $\ung=sl_2$.

\proclaim{Definition 3.1} Set $P_0=1$ and define $P_n\in{U}_q(\hat{sl}_2)^0$ by induction on $n\in\Bbb N$ using the formula
$$P_n=\frac{-k^{-1}}{1-q^{-2n}}\sum_{r=1}^{n}\psi_{r}^+P_{n-r}.$$
Define $P_{-n}=\Omega(P_n)$.\endproclaim

\vskip6pt

Definition 3.1 can be conveniently reformulated by introducing the following elements of the algebra ${U}_q(\hat{sl}_2)[[u]]$ of formal power series in an indeterminate $u$ with coefficients in ${U}_q(\hat{sl}_2)$:
$$\Psi^\pm(u)=\sum_{n=0}^\infty \psi_{\pm n}^\pm u^n,\ \ \ \ \ 
{\Cal P}^\pm(u)=\sum_{n=0}^\infty P_{\pm n}u^n.$$
Then, 3.1 is equivalent to 
$$\Psi^\pm(u)=k^{\pm 1}\frac{{\Cal P}^\pm(q^{\mp 2}u)}{{\Cal P}^\pm(u)},\tag7$$
as one sees in the case of the upper sign by multiplying both sides of (7) by ${\Cal P}^+(u)$ and equating coefficients of $u^n$ (and in the case of the lower sign by applying $\Omega$). 

\proclaim{Lemma 3.2} For all $n\in\Bbb N$, we have
$$P_n=-\frac1n \sum_{r=1}^n \frac{rq^r}{[r]_q} h_rP_{n-r}.$$
\endproclaim
\demo{Proof} We first prove that
$${\Cal P}^+(u)=\prod_{n=1}^\infty 
{\roman{exp}}\left(-\frac{q^nh_nu^n}{[n]_q}\right).\tag8$$
For this, it suffices to prove that if we substitute this formula for ${\Cal P}^+$ into the right-hand side of (7), then (7) is satisfied:
$$\align
k\frac{{\Cal P}^+(q^{-2}u)}{{\Cal P}^+(u)}&=
k\frac{\prod_{n=1}^\infty 
{\roman{exp}}\left(-\frac{q^{-n}h_nu^n}{[n]_q}\right)}
{\prod_{n=1}^\infty 
{\roman{exp}}\left(-\frac{q^nh_nu^n}{[n]_q}\right)}\\
&=k\prod_{n=1}^\infty 
{\roman{exp}}\left(\frac{(q^n-q^{-n})h_nu^n}{[n]_q}\right)\\
&=k\,{\roman{exp}}\left((q-q^{-1})\sum_{n=1}^\infty h_nu^n\right),
\endalign$$
which equals $\Psi^+(u)$ by 1.2.

Differentiating both sides of the equation 
$$\log {\Cal P}^+(u)=-\sum_{n=1}^\infty\frac{q^nh_nu^n}{[n]_q}$$
with respect to $u$ now gives
$$u\frac{({\Cal P}^+)'(u)}{{\Cal P}^+(u)}=
-\sum_{n=1}^\infty\frac{nq^nh_nu^n}{[n]_q}.$$
Multiplying both sides of this equation by ${\Cal P}^+(u)$ and then equating coefficients of $u^n$ gives the identity in 3.2. \qed\enddemo

\vskip6pt\noindent{\it Remark} The identity (8) is equivalent to 
$$P_n=q^n\sum_{{k_1,k_2,k_3\ldots\ge 0}\atop{\sum rk_r=n}}
\frac1{k_1!k_2!k_3!\cdots}
\left(-\frac{h_1}{[1]_q}\right)^{k_1}\left(-\frac{h_2}{[2]_q}\right)^{k_2}
\left(-\frac{h_3}{[3]_q}\right)^{k_3}\cdots\tag9$$
for $n\ge 0$. 
\vskip12pt We now study certain commutation formulas between the $P_n$ and the $x_r^\pm$. 

\proclaim{Lemma 3.3} Let $n\in\Bbb N$, $r\in\Bbb Z$. We have
$$P_nx_r^+=x_r^+P_n-(q^2+1)x_{r+1}^+P_{n-1}+q^2x_{r+2}^+P_{n-2}$$
(if $n=1$, the last term on the right-hand side is omitted).
\endproclaim

\vskip6pt\noindent{\it Remark} One can obtain a similar formula for $n\le 0$ by applying $\Phi\Omega$ to both sides.
\vskip6pt
\demo{Proof of 3.3} We proceed by induction on $n$. If $n=0$ there is nothing to prove, and if $n=1$ the identity follows from 
$$P_1=-qh_1,\ \ \ [h_1,x_r^+]=[2]_qx_{r+1}^+.$$

Assume then that $n\ge 2$ and that the result is known for smaller values of $n$. Then we have, by 3.2,
$$\align
[P_n,x_r^+]&=-\frac1n \sum_{m=1}^{n-1} \frac{mq^m}{[m]_q} [h_mP_{n-m},x_r^+]\\
&=-\frac1n \sum_{m=1}^{n-1}\frac{mq^m}{[m]_q}\frac{[2m]}{m}x_{r+m}^+P_{n-m}\\
&\ \ \ -\frac1n\sum_{m=1}^{n-1}\frac{mq^m}{[m]_q}h_m
\left\{-(q^2+1)x_{r+1}^+P_{n-m-1}+q^2x_{r+2}^+P_{n-m-2}\right\}
\endalign$$
by the induction hypothesis. So
$$\align
[P_n,x_r^+]&=-\frac1n\sum_{m=1}^{n-1}q^m\frac{[2m]_q}{[m]_q}x_{r+m}^+P_{n-m}\\
&\ \ \ +\frac1n\sum_{m=1}^{n-1}q^m\frac{[2m]_q}{[m]_q}
\left\{(q^2+1)x_{r+m+1}^+P_{n-m-1}-q^2x_{r+m+2}^+P_{n-m-2}\right\}\\
&\ \ \ +\frac1n (q^2+1)\sum_{m=1}^{n-1}\frac{mq^m}{[m]_q}x_{r+1}^+h_mP_{n-m-1}\\
&\ \ \ -\frac1n\sum_{m=1}^{n-1}\frac{mq^{m+2}}{[m]_q}x_{r+2}^+h_mP_{n-m-2}.
\endalign$$
If $m\ge 3$, the expression which right-multiplies $x_{r+m}^+$ on the right-hand side of the last equation is
$$-\frac1n\left\{q^m\frac{[2m]_q}{[m]_q}-q^{m-1}(q^2+1)\frac{[2m-2]_q}{[m-1]_q}+q^m\frac{[2m-4]_q}{[m-2]_q}\right\}P_{n-m}=0,$$
while the expression which right-multiplies $x_{r+1}^+$ is, by 3.2, 
$$\align
-\frac1n (q^2+1)P_{n-1} & +\frac1n (q^2+1)\sum_{m=1}^{n-1}\frac{mq^m}{[m]_q}h_mP_{n-m-1}\\
&=-\frac1n(q^2+1)P_{n-1}-\frac{(n-1)}n(q^2+1)P_{n-1}
=-(q^2+1)P_{n-1},\endalign$$
and that which right-multiplies $x_{r+2}^+$ is
$$\align
-\frac1n q^2\frac{[4]_q}{[2]_q}P_{n-2} & -\frac1n \sum_{m=1}^{n-2}\frac{mq^{m+2}}{[m]_q}h_mP_{n-m-2}
+\frac1n(q^2+1)^2P_{n-2}\\
&=\frac2n q^2P_{n-2}+\frac{(n-2)}n q^2P_{n-2}=q^2P_{n-2}.\endalign$$
This completes the inductive step, and the proof of 3.3.
\qed\enddemo

\proclaim{Lemma 3.4} Let $r\in\Bbb Z$, $n\in\Bbb N$. We have
$$x_r^+P_n=\sum_{m=0}^n q^m[m+1]_qP_{n-m}x_{r+m}^+.$$
\endproclaim

The proof is similar to that of 3.3. We omit the details.

We shall also need the following result, which is easily deduced from the $r=0$ case of 3.4:

\proclaim{Lemma 3.5} Let $n,r\in\Bbb N$. Then,
$$\!\!\!\!\!\!\!(x_0^+)^{r}P_n=\sum q^{m_1+m_2+\cdots+m_r}[m_1+1]_q[m_2+1]_q\ldots[m_r+1]_qP_{n-m_1-m_2-\cdots-m_r}x_{m_1}^+x_{m_2}^+\ldots x_{m_r}^+,$$
where the sum is over those non-negative integers $m_1,m_2,\ldots,m_r$ such that 

\noindent $m_1+m_2+\cdots+m_r\le n$.
\endproclaim
Let $U_q^{\roman{res}}(\hat\ung)^\pm$ be the $\Bbb C[q,q^{-1}]$-subalgebras of $\uqgh$ generated by the $(x_{i,r}^\pm)^{(n)}$ for all $i\in I$, $r\in\Bbb Z$, $n\in\Bbb N$, and let $U_q^{\roman{res}}(\hat\ung)^0$ be the $\Bbb C[q,q^{-1}]$-subalgebra of $\uqgh$ generated by $\uqresh^0$ and the $P_{i,r}=\varphi_i(P_r)$ for all $i\in I$, $r\in\Bbb Z$. Note that, by (9), $U_q^{\roman{res}}(\hat\ung)^0$ is also generated by $\uqres^0$ and the $h_{i,r}/[r]_{q_i}$ for all $i\in I$, $r\in\Bbb Z\backslash\{0\}$. Since, by 1.2, $(x_{i,r}^\pm)^{(n)}$ is equal, up to a sign, to a power of $T_{\check{\lambda}_i}$ applied to $(e_i^\pm)^{(n)}$, it is obvious that $U_q^{\roman{res}}(\hat\ung)^\pm\subset \uqresh$. That $U_q^{\roman{res}}(\hat\ung)^0\subset\uqresh$ will be proved in Section 5 (see 5.2).

\vskip24pt\noindent{\bf 4. More Identities}
\vskip12pt\noindent Let $\xi$ be an indeterminate, and form the polynomial algebra $\Bbb C(q)[\xi]$ over $\Bbb C(q)$. Define the following elements of the algebra $U_q(\hat{sl}_2)[[u]]$:
$${\Cal X}^+(u)=\sum_{n=0}^\infty x_n^+u^n,\ \ \ 
{\Cal X}^-(u)=\sum_{n=0}^\infty x_{n+1}^-u^n.$$
Let 
$${\Cal D}^\pm(u):\Bbb C(q)[\xi]\to U_q(\hat{sl}_2)[[u]]$$
be the unique homomorphisms of $\Bbb C(q)$-algebras such that
$${\Cal D}^\pm(u)(\xi)={\Cal X}^\pm(u).$$
Writing 
$${\Cal D}^\pm(u)=\sum_{n=0}^\infty D_n^\pm u^n,$$
the homomorphism property of ${\Cal D}^\pm(u)$ is equivalent to
$$D_n^\pm(fg)=\sum_{m=0}^n D_m^\pm(f)D_{n-m}^\pm(g)$$
for all $f,g\in\Bbb C(q)[\xi]$. The $\Bbb C(q)$-linear maps $D_n^\pm:\Bbb C(q)[\xi]\to U_q(\hat{sl}_2)^\pm$ are uniquely determined by this relation, together with
$$D_n^+(\xi)=x_n^+,\ \ \ D_n^-(\xi)=x_{n+1}^-, \ \ \ D_n^\pm(1)=\delta_{n,0}.$$

\proclaim{Proposition 4.1} Let $r,n\in\Bbb N$. Writing $\xi^{(r)}$ for $\xi^r/[r]_q!$, we have:
\vskip6pt\noindent(a)  $$q^{n+r-1}[n]_qD_{n}^+(\xi^{(r)})=\sum_{t=0}^n q^t[t]_qx_t^+D_{n-t}^+(\xi^{(r-1)}).$$
\vskip 6pt
\noindent (b) $$[n+r]_qD_{n}^-(\xi^{(r)})=\sum_{t=0}^n q^{-t}[n-t+1]_qD_t^-(\xi^{(r-1)})x_{n-t+1}^-.$$
\endproclaim

Let $T$ be the $\Bbb C(q)$-algebra automorphism $T_{\check{\lambda}_1}$ of $U_q(\hat{sl}_2)$ defined in Section 1; we have 
$$T(x_n^\pm)=-x_{n\mp 1}^\pm,\ \ \ T(\psi_n^\pm)=\psi_n^\pm,\ \ \  T(c^{1/2})=c^{1/2}$$
for all $n\in\Bbb N$.

\proclaim{Proposition 4.2} Let $r,n\ge 1$. Then,
$$D_n^+(\xi^{(r)})=\sum_{s=1}^{r-1}(-1)^{s+1}q^{s(r-1)}(x_0^+)^{(s)}
D_n^+(\xi^{(r-s)})+q^{r(r-1)}TD_{n-r}^+(\xi^{(r)}),$$
the last term being present if and only if $r\le n$.\endproclaim

The main tool for proving these identities is the next lemma. Let $U_q(\hat{sl}_2)^{++}$ be the $\Bbb C(q)$-subalgebra of $U_q(\hat{sl}_2)$ generated by the $x_n^+$ for all $n\ge 0$. It follows from the Poincar\'e--Birkhoff--Witt basis of $U_q(\hat{sl}_2)$ given in [1], Proposition 6.1, that 
$$\{(x_n^+)^{s_n}(x_{n-1}^+)^{s_{n-1}}\ldots (x_1^+)^{s_1}(x_0^+)^{s_0}\}_{n\ge 0,s_0,s_1,\ldots,s_{n-1},s_n\ge 0}$$
is a $\Bbb C(q)$-basis of $U_q(\hat{sl}_2)^{++}$. An element $x\in U_q(\hat{sl}_2)^{++}$ is said to have {\it degree} $r\ge 1$ if $x$ is a linear combination of monomials $(x_n^+)^{s_n}(x_{n-1}^+)^{s_{n-1}}\ldots (x_1^+)^{s_1}(x_0^+)^{s_0}$ with $s_0+\ldots+s_n=r$. 

\proclaim{Lemma 4.3} Let $r\ge 1$ and let $x\in U_q(\hat{sl}_2)^{++}$ have degree $r$. Suppose that $x$ acts as zero on $v_1^{\otimes r}\in V(1)_{a_1}\otimes V(1)_{a_2}\otimes\cdots\otimes V(1)_{a_r}$ for all non-zero $a_1,a_2,\ldots,a_r\in \Bbb C(q)$. Then, $x=0$.\endproclaim
\demo{Proof of 4.3} We proceed by induction on $r$. If $r=1$, then
$$x=\sum_{n=0}^\infty \lambda_nx_n^+,$$
with all but finitely many coefficients $\lambda_n\in\Bbb C(q)$ being zero. Using the formulas in Section 2 we have, in $V(1)_{a}$,
$$x.v_1=\left(\sum_{n=0}^\infty\lambda_na^n\right)v_0.$$
This can vanish for all non-zero $a\in\Bbb C(q)$ only if $\lambda_n=0$ for all $n$. 

Now assume the result for $x\in U_q(\hat{sl}_2)^{++}$ of degree $<r$, and consider
$$x=\sum \lambda(s_0,s_1,\ldots,s_n)(x_n^+)^{s_n}\ldots(x_0^+)^{s_0},\tag10$$
where the sum is over those $s_0,s_1,\ldots\ge 0$ such that $\sum_ms_m=r$, and all but finitely many $\lambda(s_0,s_1,\ldots)\in\Bbb C(q)$ are zero. Let $N$ be the maximum value of $n$ that occurs in (10), i.e.
$$N={\roman{max}}\{n\,\mid\, \text{there exists $s_0,\ldots,s_n$ with $\lambda(s_0,\ldots,s_n)\ne 0$}\},$$
and consider the action of $x$ on $V(1)_{a_1}\otimes V$, where $V=V(1)_{a_2}\otimes\cdots\otimes V(1)_{a_r}$. Clearly, 
$$x.v_1^{\otimes r}=F(a_1,a_2,\ldots,a_r)v_0^{\otimes r},$$
where $F(a_1,a_2,\ldots,a_r)$ is a polynomial in $a_1,a_2,\ldots,a_r$ with coefficients in $\Bbb C(q)$ whose degree in $a_1$ is $\le N$. We are going to identify the terms in $F(a_1,a_2,\ldots,a_r)$ that involve $a_1^N$. Since $x$ acts as zero on $v_1^{\otimes r}$, the coefficient of $a_1^N$ must actually be zero.

Let $X^\pm$ be the linear subspace of $U_q(\hat{sl}_2)$ spanned (over $\Bbb C(q)$) by the $x_n^\pm$ for all $n\in\Bbb Z$. By Proposition 5.4 in [4], if $n\ge 0$ we have
$$\Delta(x_n^+)\equiv \sum_{m=0}^n x_m^+\otimes\psi_{n-m}^+ +1\otimes x_n^+\ \ \text{modulo $U_q(\hat{sl}_2)(X^+)^2\otimes U_q(\hat{sl}_2)X^-$}.$$
Since $(X^+)^2$ annihilates $V(1)_{a_1}$ for all $a_1$, it follows that $x_n^+$ acts on $V(1)_{a_1}\otimes V$ as $x_n^+\otimes k+1\otimes x_n^+$ and hence, by an easy induction on $s\ge 1$, that $(x_n^+)^s$ acts as
$$(x_n^+\otimes k+1\otimes x_n^+)^s=\sum_{t=0}^s q^{t(s-t)}\left[s\atop t\right]_q(x_n^+)^t\otimes(x_n^+)^{s-t}k^t.\tag11$$
It is now clear that terms involving $a_1^N$ in $x.(v_1\otimes v_1^{\otimes r-1})$ can arise only from the part
$$\sum_{s_0,\ldots,s_N}\lambda(s_0,\ldots,s_N)\Delta(x_N^+)^{s_N}\Delta((x_{N-1}^+)^{s_{N-1}}\ldots (x_0^+)^{s_0})\tag12$$
of $\Delta(x)$, and in (12) only from the terms in which $x_N^+$ occurs in the first factor in the tensor product (note that if a square or higher power of $x_N^+$ occurs in the first factor of the tensor product, the corresponding term acts as zero on $v_1\otimes v_1^{\otimes r-1}$). By (11), the coefficient of $a_1^N$ in $x.v_1^{\otimes r}$ is therefore
$$\align
\sum_{s_0,\ldots,s_N}&\lambda(s_0,\ldots,s_N)q^{s_N-1}[s_N]_qq^{r-2s_N+1}
v_1\otimes(x_N^+)^{s_N-1}(x_{N-1}^+)^{s_{N-1}}\ldots(x_0^+)^{s_0}
.v_1^{\otimes r-1}\\
&=
q^rv_1\otimes\sum_{s_0,\ldots,s_N}q^{-s_N}[s_N]_q\lambda(s_0,\ldots,s_N)(x_N^+)^{s_N-1}(x_{N-1}^+)^{s_{N-1}}\ldots(x_0^+)^{s_0}.v_1^{\otimes r-1}.
\endalign$$
Thus, if $x$ acts as zero on $v_1^{\otimes r}$, this last expression must vanish, so by the induction hypothesis,
$$\sum_{s_0,\ldots,s_N}q^{-s_N}[s_N]_q\lambda(s_0,\ldots,s_N)(x_N^+)^{s_N-1}(x_{N-1}^+)^{s_{N-1}}\ldots(x_0^+)^{s_0}=0.$$
By [1], Proposition 6.1, this forces all the coefficients $\lambda(s_0,\ldots,s_N)$ to vanish. But this contradicts the assumption that $x_N$ does occur on the right-hand side of (10).\qed\enddemo

We shall also need the computations contained in the next two lemmas. 

\proclaim{Lemma 4.4} Let $r\ge s\ge 1$, and let $a_1,\ldots,a_r\in\Bbb C(q)$ be non-zero. Then, in $V(1)_{a_1}\otimes\cdots\ot V(1)_{a_r}$, we have 
$$(x_0^+)^{(s)}({\Cal X}^+(u))^{(r-s)}.v_1^{\otimes r}=f_{s,r}(a_1,\ldots,a_r;u)v_0^{\otimes r},\tag13$$
where 
$$f_{s,r}(a_1,\ldots,a_r;u)=\frac1{(1-a_1u)\cdots(1-a_ru)}\sum_{t=0}^s(-1)^t
\left[{r-t}\atop{s-t}\right]_qq^{t(r-s)}\Cal E_t(a_1,\ldots,a_r)u^t$$
and $\Cal E_t(a_1,\ldots,a_r)$ is the $t^{th}$ elementary symmetric function of $a_1,\ldots,a_r$, i.e.
$$\Cal E_t(a_1,\ldots,a_r)=\sum_{1\le r_1<r_2<\cdots<r_t\le r}
a_{r_1}a_{r_2}\ldots a_{r_t}.$$
\endproclaim
\demo{Proof} It is clear that (13) holds with some scalars $f_{s,r}(a_1,\ldots,a_r;u)$ that are formal power series in $u$ with coefficients in $\Bbb C(q)$. We can derive an inductive formula for them as follows.

Working in $V(1)_{a_1}\otimes V(1)_{a_2}\otimes\cdots\otimes V(1)_{a_r}$ we have, by (11) and the fact that $V(1)_{a_1}$ is annihilated by $(X^+)^2$,
$$(x_0^+)^{(s)}{\Cal X}^+(u)^{r-s}.v_1^{\otimes r}=(1\otimes (x_0^+)^{(s)}+q^{s-1}x_0^+\otimes(x_0^+)^{(s-1)}k){\Cal X}^+(u)^{r-s}.(v_1\otimes v_1^{\otimes r-1}).$$
Arguing as in the proof of 4.3, we see that ${\Cal X}^+(u)$ acts on $V(1)_{a_1}\otimes V$ as 

\noindent ${\Cal X}^+(u)\otimes\Psi^+(u)+1\otimes {\Cal X}^+(u)$. Hence,
$$\align
(x_0^+)^{(s)}{\Cal X}^+(u)^{r-s}.v_1^{\otimes r}&={\Cal X}^+(u).v_1\otimes (x_0^+)^{(s)}\left(\sum_{t=0}^{r-s}{\Cal X}^+(u)^t\Psi^+(u){\Cal X}^+(u)^{r-s-t-1}\right).v_1^{\otimes r-1}\\
&\ \ \ \ \ \ \ \ \ \ \ \ +q^{-s+1}v_0\otimes k(x_0^+)^{(s-1)}{\Cal X}^+(u)^{r-s}
.v_1^{\otimes r-1}\\
&=A+B,\tag14\endalign$$
say. Clearly,
$$B=q^{r-s}f_{s-1,r-1}(a_2,\ldots,a_r;u)v_0^{\otimes r}.$$

To evaluate $A$, we first show that, for $r\ge 1$, 
$$\align
\sum_{t=0}^r{\Cal X}^+(u)^t\Psi^+(u){\Cal X}^+(u)^{r-t}&\equiv q^r[r+1]_qk^{-1}{\Cal X}^+(u)^r-(q-q^{-1})x_0^-{\Cal X}^+(u)^{r+1}\\
&\ \ \ \ \ \ \ \ \ \ \ \ \ \ \ \ \ \ \ \ \ \ \ \ \ \ \text{modulo $U_q(\hat{sl}_2)X^-[[u]]$}.\tag15\endalign$$
The proof is by induction on $r$. If $r=1$, we have
$$[{\Cal X}^+(u),x_0^-]=\frac{\Psi^+(u)-k^{-1}}{q-q^{-1}},$$
so
$$\align
\Psi^+(u){\Cal X}^+(u)&=(q-q^{-1}){\Cal X}^+(u)x_0^-{\Cal X}^+(u)-(q-q^{-1})x_0^-{\Cal X}^+(u)^2+k^{-1}{\Cal X}^+(u)\\
&=(q-q^{-1}){\Cal X}^+(u)\left({\Cal X}^+(u)x_0^- -\frac{\Psi^+(u)-k^{-1}}{q-q^{-1}}\right)\\
&\ \ \ \ \ \ \ \ \ \ \ \ \ \ \ \ \ \ \ \ \ \ \ \ \ \ -(q-q^{-1})x_0^-{\Cal X}^+(u)^2+k^{-1}{\Cal X}^+(u).\endalign$$
Hence,
$$\Psi^+(u){\Cal X}^+(u)+{\Cal X}^+(u)\Psi^+(u)=q[2]_qk^{-1}{\Cal X}^+(u)-(q-q^{-1})x_0^-{\Cal X}^+(u)^2+(q-q^{-1}){\Cal X}^+(u)^2x_0^-,$$
proving (15) when $r=1$. Assume now that (15) is known for some $r\ge 1$. Then,
$$\align
\sum_{t=0}^{r+1}&{\Cal X}^+(u)^t\Psi^+(u){\Cal X}^+(u)^{r-t+1}\\
&={\Cal X}^+(u)\sum_{t=0}^r{\Cal X}^+(u)^t\Psi^+(u){\Cal X}^+(u)^{r-t}+\Psi^+(u){\Cal X}^+(u)^{r+1}\\
&\equiv {\Cal X}^+(u)(q^r[r+1]_qk^{-1}{\Cal X}^+(u)^r-(q-q^{-1})x_0^-{\Cal X}^+(u)^{r+1})
+\Psi^+(u){\Cal X}^+(u)^{r+1}\\
&\ \ \ \ \ \ \ \ \ \ \ \ \ \ \ \ \ \ \ \ \ \ \ \ \ \ \ \ \ \text{modulo $U_q(\hat{sl}_2)X^-[[u]]$}\\
&\equiv q^{r+2}[r+1]_qk^{-1}{\Cal X}^+(u)^{r+1}-((q-q^{-1})x_0^-{\Cal X}^+(u)+\Psi^+(u)-k^{-1}){\Cal X}^+(u)^{r+1}\\
&\ \ \ \ \ \ \ \ \ \ \ \ \ \ \ \ +\Psi^+(u){\Cal X}^+(u)^{r+1}\ \ \text{modulo $U_q(\hat{sl}_2)X^-[[u]]$}\\
&\equiv q^{r+1}[r+2]_qk^{-1}{\Cal X}^+(u)^{r+1}-(q-q^{-1})x_0^-{\Cal X}^+(u)^{r+2}\ \ \text{modulo $U_q(\hat{sl}_2)X^-[[u]]$},\endalign$$
completing the inductive step. 

Thus,
$$A={\Cal X}^+(u).v_1\otimes(x_0^+)^{(s)}\left(q^{r-s-1}[r-s]_qk^{-1}{\Cal X}^+(u)^{r-s-1}-(q-q^{-1})x_0^-{\Cal X}^+(u)^{r-s}\right).v_1^{\otimes r-1}.$$
Now, by the formulas in 2.3, we have, in $V(1)_{a_1}$,
$${\Cal X}^+(u).v_1=\sum_{n=0}^\infty u^nx_n^+.v_1=\sum_{n=0}^\infty a_1^nu^nv_0=\frac1{1-a_1u}v_0.$$
Moreover, an easy induction on $s$ shows that 
$$(q-q^{-1})[(x_0^+)^{(s)},x_0^-]=(q^{-s+1}k-q^{s-1}k^{-1})(x_0^+)^{(s-1)}.$$
Using these results, we get
$$\align
(1-a_1u)A&=v_0\otimes\left(q^{r-s-1+2(r-s-1)-2(r-1)}[r-s]_q(x_0^+)^{(s)}{\Cal X}^+(u)^{r-s-1}\right.\\
&\ \ \ \ \ \ \ \ \ \ \ \ \ \ \ \ \ \ \ \left.-(q^{-s+1+r-1}-q^{s-1-r+1})(x_0^+)^{(s-1)}{\Cal X}^+(u)^{r-s-1}\right).v_1^{\otimes r-1}\\
&=\left(q^s[r-s]_qf_{s,r-1}(a_2,\ldots,a_r;u)-(q^{r-s}-q^{-r+s})f_{s-1,r-1}(a_2,\ldots,a_r;u)\right)v_0^{\otimes r}.\endalign$$

Inserting these values of $A$ and $B$ into (14) gives
$$\aligned
f_{s,r}(a_1,\ldots,a_r;u)&=\left(\frac{q^{-r+s}-q^{r-s}a_1u}{1-a_1u}\right)f_{s-1,r-1}(a_2,\ldots,a_r;u)\\
&\ \ \ \ \ \ \ \ \ \ \ \ \ +q^s[r-s]_qf_{s,r-1}(a_2,\ldots,a_r;u).\endaligned\tag16$$
This identity clearly determines the $f_{s,r}$, by induction on $r$, in terms of $f_{0,1}$ and $f_{1,1}$. It therefore suffices to prove that the formula for $f_{s,r}$ in the statement of the lemma satisfies (16) and is correct when $r=1$.

Correctness for $r=1$ is trivially checked. To verify (16), we must show that
$$\align
\sum_{t=0}^s(-1)^tq^{t(r-s)}&\left[{r-s}\atop{t-s}\right]_q\Cal E_t(a_1,\ldots,a_r)u^t -q^s\sum_{t=0}^s(-1)^tq^{t(r-s-1)}\Cal E_t(a_2,\ldots,a_r)u^t\\
&=(q^{-r+s}-q^{r-s}a_1u)\sum_{t=0}^{s-1}(-1)^tq^{t(r-s)}\left[{r-t-1}\atop{s-t-1}\right]_q\Cal E_t(a_2,\ldots,a_r)u^t.\endalign$$
Thus, equating coefficients of $u^t$, we are reduced to proving that
$$\align
\left[{r-t}\atop{s-t}\right]_q\Cal E_t(a_1,\ldots,a_r)&-q^{s-t}\left[{r-t-1}\atop{s-t}\right]_q\Cal E_t(a_2,\ldots,a_r)\\
&=q^{-r+s}\left[{r-t-1}\atop{s-t-1}\right]_q\Cal E_t(a_2,\ldots,a_r)+\left[{r-t}\atop{s-t}\right]_qa_1\Cal E_{t-1}(a_2,\ldots,a_r).\endalign$$
Since 
$$\Cal E_t(a_1,\ldots,a_r)=a_1\Cal E_{t-1}(a_2,\ldots,a_r)+\Cal E_t(a_2,\ldots,a_r),$$
this reduces to
$$\left[{r-t}\atop{s-t}\right]_q-q^{s-t}\left[{r-t-1}\atop{s-t}\right]_q
=q^{-r+s}\left[{r-t-1}\atop{s-t-1}\right]_q.$$
This identity is easily checked.\qed\enddemo

\proclaim{Lemma 4.5} Let $r\ge 1$ and let $a_1,\ldots,a_r\in\Bbb C(q)$ be non-zero. Then, in 

\noindent $V(1)_{a_1}\otimes\cdots\otimes V(1)_{a_r}$, we have
$$T({\Cal D}^+(u)(\xi^{(r)})).v_1^{\otimes r}=a_1\ldots a_r{\Cal D}^+(\xi^{(r)}).v_1^{\otimes r}.$$\endproclaim
\demo{Proof} Assume that, for all $1\le s\ne t\le r$, $a_s/a_t\notin q^{\Bbb Z}$. Then, \hbox{$V=V(1)_{a_1}\otimes\cdots\otimes V(1)_{a_r}$} is an irreducible highest weight representation of $U_q(\hat{sl}_2)$ (see [9]). Let $\rho:U_q(\hat{sl}_2)\to{\roman{End}}(V)$ denote the action. Then, $\rho\circ T$ is another irreducible representation of $U_q(\hat{sl}_2)$, and since $T$ takes $U_q(\hat{sl}_2)^+$ to itself and fixes $U_q(\hat{sl}_2)^0$ pointwise, $\rho\circ T$ is highest weight with the same highest weight as $\rho$. Hence, $\rho\circ T$ is equivalent to $\rho$, i.e. there exists $F\in{\roman{GL}}(V)$ such that
$$F(x.v)=T(x).F(v)\ \ \ \text{for all $x\in U_q(\hat{sl}_2)$, $v\in V$.}$$
We might as well assume that $F(v_0^{\otimes r})=v_0^{\otimes r}$. Then, 
$$F(v_1^{\otimes r})=cv_1^{\otimes r}$$
for some non-zero $c\in\Bbb C(q)$. We proceed to compute $c$.

Using (11), we have
$$\align
(x_0^+)^r.v_1^{\otimes r}&=(x_0^+\otimes k+1\otimes x_0^+)^r.v_1^{\otimes r}\\
&=q^{r-1}[r]_q(x_0^+\otimes(x_0^+)^{r-1}k).(v_1\otimes v_1^{\otimes r-1})\\
&=[r]_qv_1\otimes(x_0^+)^{r-1}.v_1^{\otimes r-1},\endalign$$
and hence by iteration,
$$(x_0^+)^{(r)}.v_1^{\otimes r}=v_0^{\otimes r}.$$
Similarly,
$$(x_0^-)^{(r)}.v_0^{\otimes r}=v_1^{\otimes r}.$$
Hence,
$$v_1^{\otimes r}=(x_0^-)^{(r)}.v_0^{\otimes r}=T(x_1^-)^{(r)}.
F(v_0^{\otimes r})=F((x_1^-)^{(r)}.v_0^{\otimes r}).$$
Obviously,
$$(x_1^-)^{(r)}.v_0^{\otimes r}=c'v_1^{\otimes r},$$
for some $c'\in\Bbb C(q)$. To compute $c'$, apply $(x_0^+)^{(r)}$ to both sides of the last equation and use Proposition 3.5 in [4]:
$$\align
c'v_0^{\otimes r}&=(-1)^rq^{-r^2}P_rk^r.v_0^{\otimes r}\\
&=(-1)^r\times\,\text{(coefficient of $u^r$ in $(1-a_1u)\ldots(1-a_ru)$)}v_0^{\otimes r}\\
&=(a_1\ldots a_r)v_0^{\otimes r},\endalign$$
where the second equality used (6). Hence, $c'=a_1\ldots a_r$ and
$$v_1^{\otimes r}=(a_1\ldots a_r)F(v_1^{\otimes r})=(a_1\ldots a_r)cv_1^{\otimes r},$$
so $c=(a_1\ldots a_r)^{-1}$. 

Finally,
$$\align
T(\Cal D^+(u)(\xi^{(r)})).v_1^{\otimes r}&=T({\Cal X}^+(u)^{(r)}).v_1^{\otimes r}\\
&=c^{-1}T({\Cal X}^+(u)^{(r)}).F(v_1^{\otimes r})\\
&=(a_1\ldots a_r)F({\Cal X}^+(u)^{(r)}.v_1^{\otimes r})\\
&=(a_1\ldots a_r)F(\Cal D^+(u)(\xi^{(r)}).v_1^{\otimes r}).\endalign$$
But $\Cal D^+(u)(\xi^{(r)}).v_1^{\otimes r}$ is obviously a scalar multiple of $v_0^{\otimes r}$, and so is fixed by $F$. Hence,
$$T(\Cal D^+(u)(\xi^{(r)})).v_1^{\otimes r}=(a_1\ldots a_r)\Cal D^+(u)(\xi^{(r)}).v_1^{\otimes r},$$
as claimed.

Although we have proved this relation only under the assumption that $a_s/a_t\notin q^{\Bbb Z}$ for all $s\ne t$, it is clear that both sides are polynomials in $a_1,\ldots,a_r$, so it must hold identically.\qed\enddemo

We are finally in a position to prove 4.1 and 4.2.

\demo{Proof of 4.1} By equating coefficients of $u^n$ on both sides, the identity in part (a) is easily seen to be equivalent to
$$[r]_q({\Cal X}^+(q^2u)-{\Cal X}^+(u)){\Cal D}^+(u)(\xi^{r-1}) =q^{r-1}({\Cal D}^+(q^2u)(\xi^r) -{\Cal D}^+(u)(\xi^r)),$$
or, since ${\Cal D}^+(u)(\xi^r) ={\Cal X}^+(u)^r$, to
$$[r]_q{\Cal X}^+(q^2u){\Cal X}^+(u)^{r-1}-q^{-1}[r-1]_q{\Cal X}^+(u)^r =q^{r-1}{\Cal X}^+(q^2u)^r.\tag17$$
Now (17) holds trivially if $r=1$, and we claim that it follows for all $r$ if it holds for $r=2$. Indeed, assuming that $r\ge 3$ and that the result is known for $r-1$, we get
$$\align q^{r-1}{\Cal X}^+(q^2u)^r&=q{\Cal X}^+(q^2u)([r-1]_q{\Cal X}^+(q^2u){\Cal X}^+(u)^{r-2}-q^{-1}[r-2]_q{\Cal X}^+(u)^{r-1})\\
&=q[r-1]_q{\Cal X}^+(q^2u)^2{\Cal X}^+(u)^{r-2}-[r-2]_q{\Cal X}^+(q^2u){\Cal X}^+(u)^{r-1}\\
&=[r-1]_q([2]_q{\Cal X}^+(q^2u){\Cal X}^+(u)-q^{-1}{\Cal X}^+(u)^2){\Cal X}^+(u)^{r-2}\\
&\ \ \ \ \ \ \ -[r-2]_q{\Cal X}^+(q^2u){\Cal X}^+(u)^{r-1}\\
&=[r]_q{\Cal X}^+(q^2u){\Cal X}^+(u)^{r-1}-q^{-1}[r-1]_q{\Cal X}^+(u)^r,\endalign $$
where the penultimate equality is from the $r=2$ case.

Thus, it suffices to prove that
$$[2]_q{\Cal X}^+(q^2u){\Cal X}^+(u)-q^{-1}{\Cal X}^+(u)^2=q{\Cal X}^+(q^2u)^2.\tag18$$
By 4.3, it is enough to prove that both sides of (18) act in the same way on $v_1\ot v_1\in V(1)_{a_1}\ot V(1)_{a_2}$, for all non-zero $a_1,a_2\in\Bbb C(q)$. Recalling that
$$\Delta({\Cal X}^+(u))={\Cal X}^+(u)\ot\Psi^+(u)+1\ot {\Cal X}^+(u)\ \ \text{modulo $(U_q(\hat{sl}_2)(X^+)^2\ot U_q(\hat{sl}_2)X^-)[[u]]$,}$$
and noting that, in $V(1)_{a}$,
$${\Cal X}^+(u).v_1=\frac1{1-au}v_0,\ \ \Psi^+(u).v_0=\frac{q-q^{-1}au}{1-au}v_0,\ \ \Psi^+(u).v_1=\frac{q^{-1}-qau}{1-au}v_1$$
by the formulas in 2.3, we have
$$\align
[2]_q{\Cal X}^+(q^2u)&{\Cal X}^+(u).(v_1\ot v_1)\\
&=[2]_q\left({\Cal X}^+(q^2u)\otimes\Psi^+(q^2u){\Cal X}^+(u)+{\Cal X}^+(u)\ot{\Cal X}^+(q^2u)\Psi^+(u)\right).(v_1\ot v_1)\\
&=[2]_q\left(\frac1{1-q^2a_1u}\frac{q-qa_2u}{1-q^2a_2u}\frac1{1-a_2u}+
\frac1{1-a_2u}\frac{q^{-1}-qa_2u}{1-a_2u}\frac1{1-q^2a_2u}\right)\,v_1\ot v_1\\
&=[2]_q\left(\frac{q}{(1-q^2a_1u)(1-q^2a_2u)}+\frac{q^{-1}}{(1-a_1u)(1-a_2u)}
\right)v_1\ot v_1,\endalign$$
$$\align
q^{-1}{\Cal X}^+(u)^2.(v_1\ot v_1)&=q^{-1}\left({\Cal X}^+(u)\ot\Psi^+(u){\Cal X}^+(u)+{\Cal X}^+(u)\ot{\Cal X}^+(u)\Psi^+(u)\right).(v_1\ot v_1)\\
&=q^{-1}\left(\frac1{1-a_1u}\frac{q-q^{-1}a_2u}{1-a_2u}\frac1{1-a_2u}+\frac1{1-a_1u}\frac1{1-a_2u}\frac{q^{-1}-qa_2u}{1-a_2u}\right)v_1\ot v_1\\
&=\frac{q^{-1}[2]_q}{(1-a_1u)(1-a_2u)}v_1\ot v_1,\endalign$$
and hence
$$q{\Cal X}^+(q^2u)^2.(v_1\ot v_1)=\frac{q[2]_q}{(1-q^2a_1u)(1-q^2a_2u)}v_1\ot v_1.$$
It is now clear that the two sides of (18) agree in their action on $v_1\ot v_1$, for all $a_1,a_2$. This proves part (a).

The identity in part (b) can be converted into an equivalent identity for $D_n^+$ by using
$$D_n^-=-T\circ\Phi\circ D_n^+,$$
and this identity can then be proved by an argument similar to that used for part (a). We omit the details.\qed\enddemo

\demo{Proof of 4.2} Multiplying both sides by $u^n$ and summing from $n=0$ to $\infty$, we see that the identity to be proved is equivalent to 
$$\!\!\!\!\!{\Cal X}^+(u)^{(r)}=\sum_{s=1}^{r-1}(-1)^{s+1}q^{s(r-1)}(x_0^+)^{(s)}\left({\Cal X}^+(u)^{(r-s)}-(x_0^+)^{(r-s)}\right)+q^{r(r-1)}u^rT({\Cal X}^+(u)^{(r)})+(x_0^+)^{(r)}.$$
As in the proof of 4.1, we compute the action of both sides on $v_1^{\ot r}\in V(1)_{a_1}\ot\cdots\ot V(1)_{a_r}$, for arbitrary non-zero $a_1,\ldots,a_r\in\Bbb C(q)$. Using 4.4 and 4.5, we see that it suffices to prove that
$$\align
\!\!\!\!\!\!\!\!\!\!\!\!\!\!\!\!\!\!\!\!\!\!\!\frac1{(1-a_1u)\ldots(1-a_ru)}&=\frac1{(1-a_1u)\ldots(1-a_ru)}\sum_{s=1}^{r-1}
\sum_{t=0}^s(-1)^{s+t+1}q^{s(r-1)+t(r-s)}\left[{r-t}\atop{s-t}\right]_q
{\Cal E}_t(a_1,\ldots,a_r)u^t\\
&\ \ \ \ \ \ \ \ +\sum_{s=1}^{r-1}(-1)^sq^{s(r-1)}\left[r\atop s\right]_q+\frac{q^{r(r-1)}a_1\ldots a_ru^r}{(1-a_1u)\ldots(1-a_ru)}+1,\endalign$$
i.e., after using the identity (1),
$$\align
0&=q^{r(r-1)}a_1\ldots a_ru^r-(-1)^rq^{r(r-1)}(1-a_1u)\ldots(1-a_ru)-1\\
&\ \ \ \ \ \ \ +\sum_{s=1}^{r-1}\sum_{t=0}^s(-1)^{s+t+1}q^{s(r-1)+t(r-s)}\left[{r-t}\atop {s-t}\right]_q\Cal E_t(a_1,\ldots,a_r)u^t.\endalign$$
The constant term on the right-hand side vanishes by identity (1) again, and the coefficient of $u^r$ obviously vanishes. If $0<t<r$, the coefficient of $u^t$ on the right-hand side is
$$\align
\left((-1)^{r+t+1}q^{r(r-1)}+\sum_{s=t}^{r-1}(-1)^{s+t+1}
\right.&\left.q^{s(r-1)+t(r-s)}
\left[{r-t}\atop {s-t}\right]_q\right)\Cal E_t(a_1,\ldots,a_r)\\
&=\sum_{p=0}^{r-t}(-1)^pq^{p(r-s-1)}\left[{r-t}\atop {p}\right]_q\Cal E_t(a_1,\ldots,a_r),\endalign$$
which vanishes by identity (1) once again.\qed\enddemo

The following corollary of 4.2 will be needed later.

\proclaim{Corollary 4.5} Let $r,n\ge 1$. Then,
$$D_n^+(\xi^{(r)})=\sum\mu(s_0,s_1,\ldots)(x_0^+)^{(s_0)}(x_1^+)^{(s_1)}\ldots,$$
where the sum is over those non-negative integers $s_0,s_1,\ldots$ such that $\sum_ts_t=r$ and $\sum_t ts_t=n$, and the coefficients $\mu(s_0,s_1,\ldots)\in\Bbb C[q,q^{-1}]$. In particular, the coefficient of $(x_n^+)^{(r)}$ in $D_{nr}^+(\xi^{(r)})$ is $q^{nr(r-1)}$. \endproclaim
\demo{Proof} The first part is immediate from 4.2. The second part follows by induction on $r$, noting that the term involving $(x_n^+)^{(r)}$ can arise only from the second term on the right-hand side of the formula for $D_{nr}^+(\xi^{(r)})$ given in 4.2.\qed\enddemo

\noindent{\it Remark} It is possible to deduce an exact formula for 
$D_n^+(\xi^{(r)})$ from 4.2, namely
$$D_n^+(\xi^{(r)})=q^{-n}\sum_\pi q^{f(\pi)+rl(\pi)}(x_0^+)^{(r-l(\pi))}(x_1^+)^{(r_1)}(x_2^+)^{(r_2)}\ldots.
\tag19$$
Here, the sum is over all sequences $\pi:r_1,r_2,\ldots$ of non-negative integers such that $r_1+2r_2+3r_3+\cdots=n$ and $r_1+r_2+\cdots\le r$. To such a sequence we associate a Young diagram which has $r_1$ rows of length $1$, $r_2$ rows of length $2$, etc., so that $l(\pi)=r_1+r_2+\cdots$ is the total number of rows, and define $f(\pi)$ to be the sum of the products of the lengths of adjacent columns. For example, if $\pi$ is the sequence $2,1,3,1,0,0,\ldots$, whose Young diagram is 
\vskip4cm\noindent
we have $l(\pi)=7$, $f(\pi)=7\times 5+5\times4+4\times 1=59$. 

To prove the formula (19), one observes that 4.2 clearly determines $D_n^+(\xi^{(r)})$, by induction on $n$ for fixed $r$, and then by induction on $r$. Thus, it suffices to prove that (19) is correct when $r=1$, which is clear, and that if we define 
$D_n^+(\xi^{(r)})$ by the right-hand side of (19), then the identity in 4.2 is satisfied. We omit the details as we shall not make use of (19) in the sequel.

\vskip12pt We conclude this section with the following result, which will be needed in \hbox{Section 6.}

\proclaim{Proposition 4.6} For all $r\in\Bbb N,n\in\Bbb Z$, we have
$$\frac1{[r+1]_q!}((x_0^+)^rx_n^++(x_0^+)^{r-1}x_n^+x_0^++\cdots
+x_n^+(x_0^+)^r)
\in U_q^{\roman{res}}(\hat{sl}_2)^+.$$
\endproclaim
\demo{Proof} For $r\in\Bbb N,n\in\Bbb Z$, define elements $A_{r,n}\in U_q(\hat{sl}_2)$ inductively by 
$$A_{0,n}=x_n^+,\ \ \ A_{r,n}=A_{r-1,n}x_0^+-q^{2r}x_0^+A_{r-1,n},\tag20$$
and let 
$$B_{r,n}=\sum_{s=0}^r(x_0^+)^{s}x_n^+(x_0^+)^{r-s}.$$
We shall prove the following identities:
$$\align
\frac{A_{r,n}}{[r+1]_q!}&=\sum_{s=0}^r(-1)^sq^{r(r-s+1)}(x_1^+)^{(s)}T
\left(\frac{B_{r-s,n-r-1}}{[r-s+1]_q!}\right),\tag21\\
\frac{B_{r,n}}{[r+1]_q!}&=\sum_{s=0}^rq^{(r-s)(s+1)}(x_0^+)^{(r-s)}
\frac{A_{s,n}}{[s+1]_q!}.\tag22\endalign$$
The statement in the proposition for $n\ge 0$ follows from these identities by induction on $n$. Indeed, the result is trivial if $n=0$, since $B_{r,0}/[r+1]_q!=(x_0^+)^{(r+1)}$. Assuming that $B_{r,m}/[r+1]_q!\in U_q^{\roman{res}}(\hat{sl}_2)^+$ for all $r\in\Bbb N$, $m<n$, it follows from (21) that $A_{r,n}/[r+1]_q!\in U_q^{\roman{res}}(\hat{sl}_2)^+$ for all $r\in\Bbb N$, and then from (22) that $B_{r,n}/[r+1]_q!\in U_q^{\roman{res}}(\hat{sl}_2)^+$.

Note that (21) can be written in the form 
$$A_{r,n}=\sum_{s=0}^r c_{r,s}(x_1^+)^sx_{n-r}^+(x_1^+)^{r-s},$$
where
$$c_{r,s}=\sum_{t=0}^s(-1)^tq^{r(r-t+1)}\left[{r+1}\atop t\right]_q.$$
Comparing coefficients of $(x_1^+)^sx_{n-r}^+(x_1^+)^{r-s}$ on both sides of the second equation in (20), we see using (3) that (21) follows from
$$c_{r,s}=q^{2r}c_{r-1,s}-c_{r-1,s-1}.$$
This is an easy consequence of the identity
$$\left[{r+1}\atop t\right]_q=q^t\left[r\atop t\right]_q+q^{t-r+1}\left[r\atop{t-1}\right]_q.$$

To prove (22), we proceed by induction on $r$. If $r=0$, the result is trivial. Assume now that (22) holds for all smaller values of $r$ and all $n$. Then,
$$\align
\frac{B_{r,n}}{[r+1]_q!}&=\frac{(x_0^+)^rx_n^++B_{r-1,n}x_0^+}{[r+1]_q!}\\
&=\frac1{[r+1]_q!}\left((x_0^+)^rx_n^++[r]_q!\left(\sum_{s=0}^{r-1}q^{(r-s-1)(s+1)}(x_0^+)^{(r-s-1)}\frac{A_{s,n}}{[s+1]_q!}\right)x_0^+\right)\\
&=\frac1{[r+1]_q!}\left((x_0^+)^rx_n^++[r]_q!\sum_{s=0}^{r-1}q^{(r-s-1)(s+1)+2s+2}x_0^+(x_0^+)^{(r-s-1)}\frac{A_{s,n}}{[s+1]_q!}\right.\\
&\ \ \ \ \ \ \ \ \ \ \ \ \ \ \ \ \ \ \left.+[r]_q!\sum_{s=0}^{r-1}q^{(r-s-1)(s+1)}(x_0^+)^{(r-s-1)}\frac{A_{s+1,n}}{[s+1]_q!}\right),\endalign$$
using the inductive relation
$$A_{s,n}x_0^+-q^{2s+2}x_0^+A_{s,n}=A_{s+1,n}.$$
Comparing coefficients of $(x_0^+)^{(r-s)}\frac{A_{s,n}}{[s+1]_q!}$, we see that to prove (22) it is enough to show that 
$$q^{(r-s)(s+1)}[r+1]_q=q^{(r-s-1)(s+1)+2s+2}[r-s]_q+q^{(r-s)s}[s+1]_q.$$
This identity is easily checked.

This completes the proof of 4.6 for $n\ge 0$. The statement for $n\le 0$ follows from this by applying $\Phi\Omega$.\qed\enddemo

\vskip36pt\noindent{\bf 5. The Basic Lemma}
\vskip12pt
\noindent In this section, we prove the following result which is central to the whole paper. We work in $U_q(\hat{sl}_2)$. 

If $x\in U_q(\hat{sl}_2)$, let $L_x:U_q(\hat{sl}_2)\to U_q(\hat{sl}_2)$ be left multiplication by $x$. Define
$$\calpd_{n}=L_{P_n}D_0^+ +L_{P_{n-1}}D_1^+ +\cdots+L_{P_0}D_n^+.$$

\proclaim{Lemma 5.1} Let $r,s\in\Bbb N$. The following identity holds in $U_q(\hat{sl}_2)$:
$$(x_0^+)^{(r)}(x_1^-)^{(s)}=\sum_{t=0}^{{\roman{min}}(r,s)}
\sum_{{m+n=t}\atop{m,n\ge 0}}(-1)^tq^{-t(r+s-t)}D_m^-(\xi^{(s-t)})k^t\calpd_n(\xi^{(r-t)}).$$
\endproclaim

\vskip6pt\noindent{\it Remarks} 1. By applying $\Omega$ and/or powers of $T$ to both sides, analogous identities can be obtained for $(x_m^+)^{(r)}(x_n^-)^{(s)}$ whenever $m+n=\pm 1$. (The case $m+n=0$ can be deduced similarly from equation (2).)

2. Lemma 5.1 is the quantum analogue of Lemma 7.5 in [9]. The classical result is much easier, however, because there the $x_n^+$ commute among themselves, as do the $x_n^-$. In addition, one has the relation
$$D_n^\pm=\frac{(D_1^\pm)^n}{n!}$$
classically, whereas in the quantum case this equation does not even make sense.

\vskip12pt

Before proving 5.1, we note the following important corollary.

\proclaim{Corollary 5.2} For all $n\in\Bbb N$, 
\vskip 6pt
$(i)_n$ $D_n^{{}\pm{}}(\xi^{(r)})\in U_q^{\roman{res}}(\hat{sl}_2)$ for all $r\in\Bbb N$;
\vskip 6pt
$(ii)_n$ $P_n\in U_q^{\roman{res}}(\hat{sl}_2)$.
\endproclaim

Of course, it follows from part (ii), 1.3 and 3.1 that  $P_n\in U_q^{\roman{res}}(\hat{sl}_2)$ for all $n<0$.

\demo{Proof} We prove both statements simultaneously by induction on $n$, according to the following scheme:
$${\roman{(i)}}^n\implies {\roman{(ii)}}^n\implies {\roman{(i)}}^{n+1},$$
where ${\roman{(i)}}^n$ is the statement that ${\roman{(i)}}_m$ holds for all $m\le n$, and similarly for ${\roman{(ii)}}^n$.

Assume that ${\roman{(i)}}^n$ holds. Applying 5.1, we see that
$$(x_0^+)^{(n)}(x_1^-)^{(n)} =(-1)^nq^{-n^2}k^nP_n +y,$$ 
where $y$ is a sum of terms already known to belong to $U_q^{\roman{res}}(\hat{sl}_2)$ by the induction hypothesis ${\roman{(i)}}^n$ and ${\roman{(ii)}}^{n-1}$.

Now assume that ${\roman{(ii)}}^n$ holds. Applying 5.1 again, we get 
$$(x_0^+)^{(n+r)}(x_1^-)^{(n)} =(-1)^nq^{-n(n+r)}k^nD_n^+(\xi^{(r)})+z,$$
where $z$ is  a sum of terms which are already known to belong to $U_q^{\roman{res}}(\hat{sl}_2)$ by the induction hypothesis ${\roman{(i)}}^n$ and ${\roman{(ii)}}^n$. Hence, $D_n^+(\xi^{(r)})\in U_q^{\roman{res}}(\hat{sl}_2)$ for all $r\in\Bbb N$. The result for $D_n^-$ follows, since
$$D_n^-=-T\circ\Phi\circ D_n^+ ,$$
and $T$ and $\Phi$ preserve $U_q^{\roman{res}}(\hat{sl}_2)$.\qed\enddemo

\vskip12pt

Lemma 5.1 is a consequence of the next lemma.
\proclaim{Lemma 5.3} Let $n,r\in\Bbb N$. The following identity holds in $U_q(\hat{sl}_2)$:
$$\calpd_n(\xi^{(r)})x_1^-=-q^{-n-r}[n+1]_qk\calpd_{n+1}(\xi^{(r-1)})+\sum_{m=1}^{n+1}q^{m-1}[m]_qx_m^-\calpd_{n-m+1}(\xi^{(r)}).$$
\endproclaim
We complete the proof of 5.1, assuming 5.3. We proceed by induction on $s$, the case $s=0$ being trivial. By the induction hypothesis and 5.3, we have 
$$\align
[s+1]_q(x_0^+)^{(r)}(x_1^-)^{(s+1)}&=\sum_{t=0}^{{\roman{min}}(r,s)}
\sum_{m+n=t\atop m,n\ge 0}(-1)^tq^{-t(r+s-t)}D_m^-(\xi^{(s-t)})k^t\calpd(\xi^{(r-t)})x_1^-\\
&=\sum_{t=0}^{{\roman{min}}(r,s)}\sum_{m+n=t\atop m,n\ge 0}(-1)^tq^{-t(r+s-t)}D_m^-(\xi^{(s-t)})k^t\\
&\ \ \ \ \ \times\,\left\{-q^{n-r+t}[n+1]_qk\calpd_{n+1}(\xi^{(r-t-1)})\right.\\
&\ \ \ \ \ \ \ \ \ \ \ \ \left.+\sum_{l=1}^{n+1}q^{l-1}[l]_qx_l^-\calpd_{n-l+1}(\xi^{(r-t)})\right\}.
\endalign$$
Looking at the expression which left-multiplies $\calpd_n(\xi^{(r-t)})$, we see that we are reduced to proving that
$$\align 
[s+1]_qD_{t-n}^-(\xi^{(s+1-t)})q^{-t(r+s+1-t)}&= D_{t-n}^-(\xi^{(s+1-t)})q^{-(t-1)(r+s+1-t)-n-r+t}[n]_q\\
&+\sum_{v=0}^{t-n}D_v^-(\xi^{(s-t)})[t-v-n+1]_qq^{-t(r+s-t)-t-v-n}x_{t-v-n+1}^-,
\endalign$$
or, on simplifying, that 
$$[s-n+1]_qD_{t-n}^-(\xi^{(s+1-t)})=\sum_{v=0}^{t-n}q^{-v}[t-v-n+1]_q
D_v^-(\xi^{(s-t)})x_{t-v-n+1}^-.$$
After suitably re-labelling the indices, one sees that this is equivalent to the identity in 4.1(b). 

To complete the proof of 5.1, we are thus reduced to giving the
\demo{Proof of 5.3} We proceed by induction on $r$. When $r=0$, the identity becomes
$$P_nx_1^-=\sum_{m=1}^{n+1}q^{m-1}[m]_qx_m^-P_{n-m+1}.$$
This follows from 2.4 (take $r=1$ in 3.4 and apply $\Phi$ to both sides, noting that $\Phi(P_n)=P_n$ for all $n\in\Bbb Z$). 

Assuming the result for $r$, we have
$$\align 
[r+1]_q\calpd_n(\xi^{(r+1)})x_1^-
&=\sum_{m=0}^nL_{P_m}D_{n-m}^+(\xi\xi^{(r)})x_1^-\\
&=\sum_{{m+l\le n}\atop {m,l\ge 0}}P_mx_l^+D_{n-m-l}^+(\xi^{(r)})x_1^-
\\
&=\sum_{{m+l\le n}\atop {m,l\ge 0}}[P_m, x_l^+]D_{n-m-l}^+(\xi^{(r)})x_1^-
+\sum_{l=0}^nx_l^+\calpd_{n-l}(\xi^{(r)})x_1^-
\\
&=\sum_{{m+l\ge n}\atop {m,l\ge 0}}
(q^2x_{l+2}^+P_{m-2}-(q^2+1)x_{l+1}^+P_{m-1})D_{n-m-l}^+(\xi^{(r)})x_{1}^- 
\\
&\ \ \ \ \ +\sum_{l=0}^nx_l^+\left\{-q^{-n-r+l}[n-l+1]_qk\calpd_{n-l+1}(\xi^{(r-1)})\right.\\
&\ \ \ \ \ \ \ \ \ \left.+\sum_{s=1}^{n-l+1}q^{s-1}[s]_qx_s^-\calpd_{n-l-s+1}(\xi^{(r)})\right\},
\tag23\endalign$$ 
on using 3.3 and the induction hypothesis. The first term on the right-hand side of (23) is equal to
$$\align
& q^2\sum_{l=0}^{n-2}x_{l+2}^+\calpd_{n-l-2}(\xi^{(r)})x_1^- 
-(q^2+1)\sum_{l=0}^{n-1}x_{l+1}^+\calpd_{n-l-1}(\xi^{(r)})x_1^-\\
&=q^2\sum_{l=0}^{n-2}x_{l+2}^+\left\{-q^{-n+l+2-r}[n-l-1]_qk\calpd_{n-l-1}
(\xi^{(r-1)})+\sum_{m=1}^{n-l-1}q^{m-1}[m]_qx_m^-\calpd_{n-l-m-1}(\xi^{(r)})
\right\}\\
&\ \ -(q^2+1)\sum_{l=0}^{n-1}x_{l+1}^+
\left\{-q^{-n+l+1-r}[n-l]_qk\calpd_{n-l}(\xi^{(r-1)})
+\sum_{m=1}^{n-l}q^{m-1}[m]_qx_m^-\calpd_{n-l-m}(\xi^{(r)})\right\},
\endalign$$
on using the induction hypothesis again. 

We now distinguish two types of term on the right-hand side of (23):
\vskip6pt\noindent(a) {\it The term involving $x_m^+k\calpd_{n-m+1}(\xi^{(r-1)})$}, which is left-multiplied by
$$-q^{-n-r+m-2}[n-m+1]_q-q^{-n-r+m}[n-m+1]_q+q^{-n+m-r-2}(q^2+1)[n-m+1]_q=0,$$
except that if $m=0$ only the first term appears, and if $m=1$ only the first and last terms appear. Hence, the contribution of this type of term is
$$q^{-n-r+1}[n]_qkx_1^+\calpd_n(\xi^{(r-1)})-q^{-n-r-2}[n+1]_qkx_0^+\calpd_{n+1}(\xi^{(r-1)}).$$
(b) {\it The term involving $x_m^+x_l^-\calpd_{n-m-l+1}(\xi^{(r)})$}, which is left-multiplied by
$$q^{l-1}[l]_q+q^{l+1}[l]_q-q^{l-1}(q^2+1)[l]_q=0,$$
except if $m=0$ or $1$, so the net contribution of this type of term is
$$q^{l-1}[l]_qx_0^+x_l^-\calpd_{n-l+1}(\xi^{(r)})-q^{l+1}[l]_qx_1^+x_l^-\calpd_{n-l}(\xi^{(r)}).$$

Hence,

$$\align
[r+1]_q\calpd_n(\xi^{(r+1)})x_1^-&=q^{-n-r+1}[n]_qkx_1^+\calpd_n(\xi^{(r-1)})-q^{-n-r-2}[n+1]_qkx_0^+\calpd_{n+1}(\xi^{(r-1)})\\
&\ \ \ +\sum_{l=1}^{n+1}q^{l-1}[l]_qx_0^+x_l^-\calpd_{n-l+1}(\xi^{(r)})
-\sum_{l=1}^n q^{l+1}[l]_qx_1^+x_l^-\calpd_{n-l}(\xi^{(r)}).
\endalign$$

On the other hand, we are trying to show that
$$\align
[r+1]_q & \calpd_n(\xi^{(r+1)})x_1^-\\
&=-q^{-n-r-1}[n+1]_q[r+1]_qk\calpd_{n+1}(\xi^{(r)})
+[r+1]_q\sum_{m=1}^{n+1}q^{m-1}[m]_qx_m^-\calpd_{n-m+1}(\xi^{(r+1)})\\
&=-q^{-n-r-1}[n+1]_q[r+1]_qk\calpd_{n+1}(\xi^{(r)})\\
&\ \ \ +\sum_{m=1}^{n+1}q^{m-1}[m]_qx_m^-
\left\{\sum_{l=0}^{n-m+1}x_l^+\calpd_{n-l-m+1}(\xi^{(r)})\right.\\
&\ \ \ \ \ \ \ \ \ 
\left.+\sum_{{l+s\le n-m+1}\atop{l,s\ge 0}}(q^2x_{s+2}^+P_{l-2}-(q^2+1)x_{s+1}^+P_{l-1})
\calpd_{n-m-l-s+1}(\xi^{(r)})\right\},\endalign$$
by repeating the argument leading to equation (23), which is equal to
$$\align
-q^{-n-r-1} & [n+1]_q[r+1]_qk\calpd_{n+1}(\xi^{(r)})+\sum_{{m+l\le n+1}\atop{l\ge 0,m\ge 1}}q^{m-1}[m]_qx_m^-x_l^+\calpd_{n-l-m+1}(\xi^{(r)})\\
&\ \ \ +\sum_{{s+l+m\le n+1}\atop{l\ge 0,m\ge 1}}q^{m-1}[m]_qx_m^-(q^2x_{s+2}^+P_{l-2}-(q^2+1)x_{s+1}^+P_{l-1})D_{n-l-m-s+1}^+(\xi^{(r)}).
\endalign$$
As above, one sees that the term $x_m^-x_l^+$, with $l\ge 0$, $m\ge 1$ and $l+m\le n+1$, survives in the sum of the last two terms only if $l=0$ or $1$, giving
$$\align
[r+1]_q & \calpd_n(\xi^{(r+1)})x_1^-\\
&=-q^{-n-r-1}[n+1]_q[r+1]_qk\calpd_{n+1}(\xi^{(r)})
+\sum_{m=1}^{n+1}q^{m-1}[m]_qx_m^-x_0^+\calpd_{n-m+1}(\xi^{(r)})\\
&\ \ \ \ \ \ \ \ -\sum_{m=1}^n q^{m+1}[m]_qx_m^-x_1^+\calpd_{n-m}(\xi^{(r)}).
\endalign$$

So we must prove that
$$\align
q^{-n-r+1}[n]_q & kx_1^+\calpd_n(\xi^{(r-1)})-
q^{-n-r-2}[n+1]_qkx_0^+\calpd_{n+1}(\xi^{(r-1)})\\
& +q^{-n-r-1}[r+1]_q[n+1]_qk\calpd_{n+1}(\xi^{(r)})\\
& +\sum_{m=1}^{n+1}q^{m-1}[m]_q\frac{\psi_m^+}{q-q^{-1}}\calpd_{n-m+1}(\xi^{(r)})\\
& -\sum_{m=1}^n q^{m+1}[m]_q\frac{\psi_{m+1}^+}{q-q^{-1}}\calpd_{n-m}(\xi^{(r)})=0.\endalign$$
The sum of the last two terms is equal to
$$\align
\sum_{m=0}^n\frac{\psi_{m+1}^+}{q-q^{-1}} & \calpd_{n-m}(\xi^{(r)})
=\sum_{{l+m\le n}\atop{l,m\ge 0}}\frac{\psi_{m+1}^+}{q-q^{-1}}P_lD_{n-l-m}^+(\xi^{(r)})\\
&=\sum_{s=0}^n\sum_{m=0}^s\frac{\psi_{m+1}^+}{q-q^{-1}}P_{s-m}D_{n-s}^+(\xi^{(r)})\\
&=-\sum_{s=0}^nq^{-s-1}[s+1]_qkP_{s+1}D_{n-s}^+(\xi^{(r)}),
\endalign$$
by 3.1. So we are reduced to proving that
$$\align
q[n]_qkx_1^+\calpd_n(\xi^{(r-1)}) & -q^{-2}[n+1]_qkx_0^+\calpd_{n+1}(\xi^{(r-1)})\\
&+q^{-1}[r+1]_q[n+1]_qk\calpd_{n+1}(\xi^{(r)})\tag24\\
&=q^{n+r}\sum_{s=0}^n q^{-s-1}[s+1]_qkP_{s+1}D_{n-s}^+(\xi^{(r)}).
\endalign$$
Now,
$$\align
x_0^+\calpd_{n+1}(\xi^{(r-1)})&=\sum_{m=0}^{n+1}x_0^+P_m
D_{n-m+1}^+(\xi^{(r-1)})\\
&=\sum_{m=0}^{n+1}\sum_{l=0}^m P_{m-l}x_l^+q^l[l+1]_qD_{n-m+1}^+(\xi^{(r-1)}),
\endalign$$
by 3.4, which is equal to
$$\sum_{m=0}^{n+1}P_m\sum_{t=0}^{n-m+1}q^t[t+1]_qx_t^+D_{n-m-t+1}^+(\xi^{(r-1)}).$$
Similarly,
$$\align x_1^+\calpd_n(\xi^{(r-1)})&=\sum_{m=0}^n\sum_{t=0}^mq^t[t+1]_qP_{m-t}x_{t+1}^+
D_{n-m}^+(\xi^{(r-1)})\\&=\sum_{m=0}^{n+1}P_m\sum_{t=0}^{n-m}q^t[t+1]_qx_{t+1}^+
\calpd_{n-t-m}^+(\xi^{(r-1)}).\endalign$$
Inserting these results in (24), and equating coefficients of $kP_m$ on both sides, we see that it suffices to prove that
$$\align 
q[n]_q\sum_{t=0}^{n-m}q^t[t+1]_qx_{t+1}^+ & D_{n-t-m}^+(\xi^{(r-1)})
+q^{-1}[r+1]_q[n+1]_qD_{n-m+1}^+(\xi^{(r)})\\
&-q^{-2}[n+1]_q\sum_{t=0}^{n-m+1}q^t[t+1]_qx_t^+D_{n-t-m+1}^+(\xi^{(r-1)})\\
&\ \ \ \ =q^{n+r-m}[m]_qD_{n-m+1}^+(\xi^{(r)}).
\endalign$$
The sum of the first and third terms is equal to
$$\align
[n]_q\sum_{t=0}^{n-m+1}q^t[t]_qx_t^+D_{n-t-m+1}^+(\xi^{(r-1)})&-q^{-2}[n+1]_q
\sum_{t=0}^{n-m+1}(1+q^{t+1}[t])x_t^+D_{n-t-m+1}^+(\xi^{(r-1)})\\
=-q^{-2}[n+1]_q[r]_qD_{n-m+1}^+(\xi^{(r)})&+\sum_{t=0}^{n-m+1}(q^t[t]_q[n]_q
-q^{t-1}[t]_q[n+1]_q)x_t^+D_{n-t-m+1}^+(\xi^{(r-1)})\\=-q^{-2}[n+1]_q[r]_qD_{n-m+1}^+
(\xi^{(r)})&-q^{-n-1}\sum_{t=0}^{n-m+1}q^t[t]_qx_t^+D_{n-t-m+1}^+(\xi^{(r-1)}).
\endalign$$
Thus, we are reduced to proving that
$$\align
(q^{-1}[r+1]_q[n+1]_q-q^{n+r-m}[m]_q & -q^{-2}[n+1]_q[r]_q)
D_{n-m+1}^+(\xi^{(r)})\\
&=q^{-n-1}\sum_{t=0}^{n-m+1}q^t[t]_qx_t^+D_{n-t-m+1}^+(\xi^{(r-1)}),
\endalign$$
i.e. that
$$\sum_{t=0}^{n-m+1}q^t[t]_qx_{t}^+D_{n-t-m+1}^+(\xi^{(r-1)})=
q^{n+r-m}[n-m+1]_qD_{n-m+1}^+(\xi^{(r)}).$$
This is the statement of 4.1(a) (with $n$ replaced by $n-m+1$).
\qed\enddemo

\vskip36pt

\noindent{\bf 6. A Triangular Decomposition}
\vskip12pt\noindent The following result is central to the construction of highest weight representations of $\uqresh$ in Sections 7 and 8.

\proclaim{Proposition 6.1} We have
$$\uqresh=\uqreshm.\uqresho.\uqreshp.$$
\endproclaim
\demo{Proof} Let $U^{\Delta}$ denote the right-hand side of the equation in 6.1. Following the strategy in [13], Section 5, define the {\it degree} of the generators of $\uqresh$ as follows:
$${\text{deg}}((x_{i,n}^{{}\pm{}})^{(r)})=r,\ \ 
{\text{deg}}(k_i^{{}\pm 1{}}) 
={\text{deg}}\left(\left[{k_i;r}\atop n\right]\right)=0,\ \ 
{\text{deg}}\left(\frac{h_{i,n}}{[n]_{q_i}}\right)=|n|.$$
Further, call a {\it monomial} a finite product of the form 
$$\!\!\!\!\!\!\!\!\!\bigl({\text{product of}} \ {(x_{i,n}^-)^{(r)}}'s\bigr)\times 
\bigl({\text{product of}}\ {k_i^{{}\pm 1}}'s,\ {\left[{k_i;r}\atop n\right]}'s, \ {\text{and}}\ \frac{h_{i,n}}{[n]_{q_i}}'s\bigr)\times  \bigl({\text{product of}} \ {(x_{i,n}^+)^{(r)}}'s\bigr),$$
and define its degree to be the sum of the degrees of the factors. Finally, any element of $U^\Delta$ that is a $\Bbb C[q,q^{-1}]$-linear combination of monomials of degree $d$ is said to have degree $d$. It follows from Proposition 6.1 in [1] that the degree is well defined. Note that, by (9), ${\roman{deg}}(P_{i,n})=|n|$. Finally, for $n\in\Bbb N$, let  $U_n^{\Delta}$ be the subspace of  $U^{\Delta}$ consisting of the elements of degree $\le n$. Proceeding by an evident induction on the degree, and using the analogue of 6.1 for $\uqres$ (see Section 1), it suffices to prove the following, for all $r,s\in\Bbb N$, $m,n\in\Bbb Z$, $i,j\in I$:

\vskip6pt\noindent (a) $\left[\frac{h_{j,n}}{[n]_{q_j}}, (x_{i,m}^+)^{(r)}\right]\in U_{r+|n|-1}^{\Delta}$;

\noindent
(b) $\left[\frac{h_{j,n}}{[n]_{q_j}}, (x_{i,m}^-)^{(r)}\right]\in U_{r+|n|-1}^{\Delta}$;

\noindent
(c) $[(x_{i,m}^+)^{(r)},  (x_{j,n}^-)^{(s)}]\in U_{r+s-1}^{\Delta}$,
\vskip6pt
\noindent where $[a,b]$ denotes $ab-ba$. 

\vskip6pt We consider only the case $n\ge 0$; the case $n\le 0$ is similar. By applying $T_{\check{\lambda}_i}^m$, it suffices to prove (a) when $m=0$. Note also that (b) follows from (a) by applying $\Omega$. To prove (a) and (b), it therefore suffices to prove that $\left[\frac{h_{j,n}}{[n]_{q_j}}, (x_{i,0}^+)^{(r)}\right]\in U^{\Delta}$ for all $i,j\in I$, $n,r>0$. Now,
$$\left[\frac{h_{j,n}}{[n]_{q_j}}\,,\,(x_{i,0}^+)^r\right]=\frac1n \frac{[na_{ji}]_{q_j}}{[n]_{q_j}}\sum_{s=1}^r(x_{i,0}^+)^{s-1}x_{i,n}^+
(x_{i,0}^+)^{r-s}.$$
Our assertion follows, since it is clear that $\frac{[na_{ji}]_{q_j}}{[n]_{q_j}}\in\Bbb C[q,q^{-1}]$, and by 4.6 we have
$$\frac1{[r]_{q_i}!}\sum_{s=1}^r(x_{i,0}^+)^{s-1}x_{i,n}^+
(x_{i,0}^+)^{r-s}\in\uqreshp.$$

To prove (c), we can of course assume that $i=j$, so we might as well work in the case $\ung=sl_2$ and drop the subscripts $i,j$ as usual. We first show that
$$(x_0^+)^{(r)}D_{n}^-(\xi^{(s)})\in  U^{\Delta}\ \ \ \text{ for all $r,s\in\Bbb N$.}\tag25$$
We proceed by induction on $n$, the statement being obvious when $n=0$. Assuming the result for $n$, we consider
$$(x_0^+)^{(r)}(x_0^+)^{(n+1)}(x_1^-)^{(n+s+1)}=
\left[{n+r+1}\atop r\right]_q(x_0^+)^{(n+r+1)}(x_1^-)^{(n+s+1)}.\tag26$$
By 5.1 and 5.2, the right-hand side of (26) is in $U^{\Delta}$, and by 5.1 the left-hand side equals
$$(x_0^+)^{(r)}D_{n+1}^-(\xi^{(s)})+w,$$
where $w$ belongs to $U^{\Delta}$ by the induction hypothesis and part (a). So $(x_0^+)^{(r)}D_{n+1}^-(\xi^{(s)})$ $\in U^{\Delta}$, and (25) is proved for $n+1$.

We now prove (c) by induction on $s$. As above, it suffices to prove that $$[(x_{m}^+)^{(r)}\,,\,(x_{n}^-)^{(s)}]\in U^{\Delta}.\tag27$$ 
We need only consider the case $m+n>0$. For, the case $m+n<0$ is similar, and the case $m+n=0$ can be reduced to the case $m=n=0$ by applying $T^m$, and that case is contained in (2). If $s=1$, then by applying $T^m$, we are reduced to proving that $(x_0^+)^{(r)}x_n^-\in U^{\Delta}$ for all $r,n>0$. This follows by taking $s=1$ in (25). Assume then that $s\ge 2$ and that (27) holds for smaller 
values of $s$. By 4.5,
$$D_{(n-1)s}^-(\xi^{(s)})=q^{(n-1)s(s-1)}(x_n^-)^{(s)}+z,$$
where $z$ is a linear combination of products
$$(x_{n_1}^-)^{(s_1)}(x_{n_2}^-)^{(s_2)}\ldots,$$
where $n_1,n_2,\ldots$ and $s_1,s_2,\ldots$ are positive integers such that 
$s_1,s_2,\ldots<s$ and $n_1s_1+n_2s_2+\cdots=ns$. Hence,
$$[(x_{m}^+)^{(r)}\,,\,(x_{n}^-)^{(s)}]=q^{-(n-1)s(s-1)}\left([(x_{m}^+)^{(r)}\,,\,D_{(n-1)s}^-(\xi^{(s)})]-[(x_{m}^+)^{(r)}\,,\,z].\right)$$
The first term on the right-hand side belongs to $U^\Delta$ by (25), and the second term belongs to $U^\Delta$ by the induction hypothesis. 

This completes the proof of 6.1.\qed\enddemo

\vskip6pt\noindent{\it Remark} The arguments above actually prove the following statement: if $r,s\in\Bbb N$, $m,n\in\Bbb Z$, $m+n\ge 0$, then $(x_m^+)^{(r)}(x_n^-)^{(s)}$ is a linear combination of products $X^-X^0X^+$, where $X^\pm\in{U}_q^{{\roman{res}}}(\hat{sl}_2)^\pm$ and $X^0$ lies in the subalgebra of ${U}_q^{{\roman{res}}}(\hat{sl}_2)^0$ generated by $k^{\pm 1}$, the $\left[{k;u}\atop v\right]$ for all $u\in\Bbb Z$, $v\in\Bbb N$, and the $P_N$ for $N\ge 0$. Analogous statements in $\uqresh$ can be obtained by applying $\varphi_i$.

\vskip24pt\noindent
{\bf 7. Representation Theory of} $U_\epsilon^{\roman{res}}$
\vskip 12pt\noindent
We first recall some facts about the representation theory of $\uepres$ (see [5] and [15]).

A representation $V$ of $\uepres$ is said to be of type I if, for all $i\in I$, $k_i$ acts semisimply on $V$ with eigenvalues in $\epsilon_i^{\Bbb Z}$ ($\epsilon_i=\epsilon^{d_i}$). Any finite-dimensional irreducible representation of $\uepres$ can be obtained from a type I representation by tensoring with a one-dimensional representation on which the $e_i^\pm$ act as zero and the $k_i$ act as $\pm 1$.

If $\lambda\in P$, the weight space $V_\lambda$ of a representation $V$ of $\uepres$ is defined by
$$V_\lambda=\left\{v\in V \bigm| k_i.v=\epsilon_i^{n_i^0}v,\ \ 
\left[k_i;0\atop \ell\right].v=n_i^1 v\right\},$$
where $n_i=\lambda(\check{\alpha}_i)$, $n=n^0+\ell n^1$, $0\le n^0<\ell$.
We have 
$$(e_i^{{}\pm{}})^{(r)}.V_\lambda\subseteq V_{\lambda\pm r\alpha_i},\ \ \ \ {\roman{dim}}(V_{w(\lambda)})={\roman{dim}}(V_\lambda),\tag28$$
for all $i\in I$, $r\in\Bbb N$, $w\in W$. 

A vector $v$ in a type I representation $V$ of $\uepres$ is said to be a highest weight vector if there exists $\lambda\in P$ such that 
$$v\in V_\lambda\ \ \text{and}\ \ (e_i^+)^{(r)}.v=0\ \ \ {\text{for all}} \  r>0, i\in I.$$
If, in addition, $V=\uepres.v$, then $V$ is said to be a highest weight representation with highest weight $\lambda$.

For any $\lambda\in P$, there exists, up to isomorphism, a unique irreducible representation $V(\lambda)$ of $\uepres$ with highest weight $\lambda$. We have
$$V(\lambda)=\bigoplus_{\mu\le\lambda}V(\lambda)_\mu.$$
Every finite-dimensional irreducible type I representation of $\uepres$ is isomorphic to $V(\lambda)$ for some (unique) $\lambda\in P^+$.

The following result of Lusztig [15] gives the structure of $V(\lambda)$ for arbitrary $\lambda\in P^+$.

\proclaim{Theorem 7.1} Let $\lambda\in P^+$.

\vskip6pt\noindent (a) If $\lambda(\check{\alpha}_i)<\ell$ for all $i\in I$, then $V(\lambda)$ is irreducible as a representation of ${U}_\epsilon^{{\roman{fin}}}(\ung)$.

\noindent (b) If $\lambda$ is divisible by $\ell$, say $\lambda=\ell\lambda^1$, then $V(\lambda)$ is isomorphic to the pull back of $\overline{V}(\lambda^1)$ by $Fr_{\epsilon}:\uepres\to\overline{U}(\ung)$. In particular, the weight space $V(\lambda)_\mu\ne 0$ only if $\lambda-\mu\in\ell Q^+$ and (hence) the $e_i^\pm$ act as zero on $V(\lambda)$ for all $i\in I$.

\noindent (c) In general, write $\lambda=\lambda^0+\ell\lambda^1$, where $0\le \lambda^0(\check{\alpha}_i)<\ell$ for all $i\in I$. Then, as representations of $\uepres$,
$$V(\lambda)\cong V(\lambda^0)\otimes V(\ell\lambda^1).$$
\endproclaim
\vskip6pt
We now give the analogous definitions for $\uepresh$. Let ${U}_\epsilon^{{\roman{res}}}(\hat\ung)^\pm$ be the subalgebras of $\uepresh$ generated by $(x_{i,n}^{{}\pm{}})^{(r)}$ for all $i\in I, n\in\Bbb Z$, $r\in\Bbb N$, and let ${U}_\epsilon^{{\roman{res}}}(\ung)^0$ be the subalgebra generated by the $ k_i^{{}\pm{1}}$, $\left[k_i;0\atop\ell\right]$ and $P_{i,n}$ for all $i\in I, n\in\Bbb Z$. 

\proclaim{Definition 7.2} A representation $V$ of $\uepresh$ is said to be of type I if $V$ is of type I as a representation of $\uepres$ and if $c^{1/2}$ acts as 1 on $V$.

If $V$ is a type I representation of $\uepresh$, a vector $v\in V$ is said to be a highest  weight vector if there exists a homomorphism of algebras $\Lambda:{U}_\epsilon^{{\roman{res}}}(\ung)^0\to\Bbb C$ such that
$$(x_{i,n}^{{}\pm{}})^{(r)}.v=0\ \ \ \ \text{for all}\  i\in I, n\in\Bbb Z, r\in\Bbb N,\ \ \text{and}\ \ x.v=\Lambda(x)v\ \ \ \ \text{for all}\  x\in {U}_\epsilon^{{\roman{res}}}(\ung)^0.$$
If, in addition, $V=\uepresh.v$, then $V$ is said to be a highest weight representation with highest weight $\Lambda$.
\endproclaim 

\vskip6pt\noindent{\it Remarks} 1. Since a $\uepresh$-highest weight vector $v\in V$ is, in particular, a $\uepres$-highest weight vector, there exists $\lambda\in P^+$ such that
$$\Lambda(k_i)=\epsilon_i^{n_i^0},\ \ \ \Lambda\left(\left[{k_i;0}\atop{\ell}\right]\right)=n_i^1,$$
where $n=\lambda(\check{\alpha}_i)$, $n_i=n_i^0+\ell n_i^1$, $0\le n_i^0<\ell$. The highest weight $\Lambda$ is determined by the $n_i$ and the complex numbers $\{\Lambda(P_{i,r})\}_{i\in I,r\in\Bbb Z}$ (with $\Lambda(P_{i,0})=1$). 

2. By the relations in 1.2 (and the fact that $c^{1/2}$ acts as $1$ on $V$), we have
$$(x_{i,n}^\pm)^{(r)}.V_\mu\subseteq V_{\mu+r\alpha_i},\ \ \ \ \uepresho.V_\mu=V_\mu\tag29$$
for all $\mu\in P$.
\vskip6pt

\proclaim{Proposition 7.3} For any algebra homomorphism $\Lambda:\uepresho\to\Bbb C$, there exists an irreducible representation $V(\Lambda)$ of $\uepresh$ with highest weight $\Lambda$. Moreover, $V(\Lambda)$ is unique up to isomorphism.\endproclaim
\demo{Proof} Let $M(\Lambda)$ be the quotient of $\uepresh$ by the left ideal generated by the $x_{i,r}^+$ for all $r\in\Bbb Z$, together with the elements $x-\Lambda(x)1$ for all $x\in\uepresho$. It is obvious that $M(\Lambda)$ is a highest weight representation of $\uepresh$, with highest weight vector $m_\Lambda$ (say) which is the image of $1$ in $M(\Lambda)$. By 6.1, $m_\Lambda$ is, up to scalar multiples, the unique vector of maximal weight for $\uepres$. It now follows by the usual arguments that $M(\Lambda)$ has a unique irreducible quotient representation $V(\Lambda)$, and that $V(\Lambda)$ is, up to isomorphism, the unique irreducible representation of $\uepresh$ with highest weight $\Lambda$.\qed\enddemo
 
\vskip6pt
The following multiplicative property of highest weights of representations of $\uepresh$ will be of crucial importance later.

\proclaim{Proposition 7.4} Let $V'$, $V''$ be type I representations of $\uepresh$ and let $v'\in V'$, $v''\in V''$ be highest weight vectors.
Then, $v'\ot v''$ is a highest weight vector in $V'\ot V''$ and
$$\calp_i^\pm(u).(v'\ot v'') =\calp_i^\pm(u).v'\otimes\calp_i^\pm(u).v''.\tag30$$
\endproclaim

For the proof, we shall need the following lemma. Recall the algebra homomorphisms $\varphi_i:U_q(\hat{sl}_2)\to\uqgh$ ($i\in I$) and define $\Delta_i$ to be the composite homomorphism
$$U_q(\hat{sl}_2) @>{\Delta} >> U_q(\hat{sl}_2)\otimes U_q(\hat{sl}_2)
@>{\varphi_i\otimes\varphi_i}>> \uqgh\otimes\uqgh,$$
where the first map is the comultiplication of $U_q(\hat{sl}_2)$. Since $\varphi_i(U_q^{{\roman{res}}}(\hat{sl}_2))\subset \uqresh$ and since $\uqresh$ is a $\Bbb C[q,q^{-1}]$-Hopf subalgebra of $\uqgh$, it follows that $$\Delta_i(U_q^{\roman{res}}(\hat{sl}_2))\subset\uqresh\otimes\uqresh.$$ 
Let
$X^{(+)}$ be the $\Bbb C[q,q^{-1}]$-span of $\{(x_{j,m}^+)^{(s)}\,\mid\,j\in I,m\in\Bbb Z,s\in\Bbb N\}$.

\proclaim{Lemma 7.5} For all $i\in I$, $r\in\Bbb N$, 
$$\Delta((x_{i,1}^-)^{(r)})-\Delta_i((x_{1}^-)^{(r)})\in\uqresh X^{(+)}\otimes\uqresh\,+\,\uqresh\otimes\uqresh X^{(+)}.$$
\endproclaim

Assuming this lemma for the moment, we give the 
\demo{Proof of 7.4} By 5.1, if $r\ge 0$ and $v=v'$, $v''$ or $v'\ot v''$,
$$(x_{i,0}^+)^{(r)}(x_{i,1}^-)^{(r)}.v = \epsilon_i^{-r^2}k_i^rP_{i,r}.v.$$
Now, using 1.1 and the isomorphism in 1.2, one finds that
$$\Delta_i(x_{0}^+)=x_{i,0}^+\ot k_i+ 1\ot x_{i,0}^+,\ \ 
\Delta_i(x_{1}^-)=x_{i,1}^-\ot 1+k_i\ot x_{i,1}^-.$$
From this, one easily deduces (by induction on $r$) that
$$\align 
\Delta_i((x_{0}^+)^{(r)})&=\sum_{s=0}^r\epsilon_i^{s(r-s)}(x_{i,0}^+)^{(s)}\ot (x_{i,0}^+)^{(r-s)}k_i^s,\\
\Delta_i((x_{1}^-)^{(r)})&=\sum_{s=0}^r\epsilon_i^{s(r-s)}k_i^{r-s}(x_{i,1}^-)^{(s)}\ot (x_{i,1}^-)^{(r-s)}.
\endalign $$
From 7.5, and the obvious fact that $\uqresh X^{(+)}\otimes\uqresh\,+\,\uqresh\otimes\uqresh X^{(+)}$, specialised to $q=\epsilon$, annihilates $v'\otimes v''$, it follows that
$$\align 
(-1)^r\epsilon_i^{-r^2} & k_i^rP_{i,r}.(v'\ot v'')=\Delta((x_{i,0}^+)^{(r)}(x_{i,1}^-)^{(r)}).(v'\otimes v'')\\
&=\sum_{s,t=0}^r(-1)^{s+t}\epsilon_i^{s(r-s)+t(r-t)}(x_{i,0}^+)^{(s)}k_i^{r-t}
(x_{i,1}^-)^{(t)}.v'\ot (x_{i,0}^+)^{r-s}k_i^s(x_{i,1}^-)^{(r-t)}.v''.\endalign$$
Now, by (29),
$$\align (x_{i,0}^+)^{(s)}(x_{i,1}^-)^{(t)}.v'&=0 \ \  {\text{if}} \ \ s>t,\\
(x_{i,0}^+)^{(r-s)}(x_{i,1}^-)^{(r-t)}.v''&=0 \ \  {\text{if}}\ \ r-s>r-t.\endalign$$
So, by 5.1 again,
$$\align 
(-1)^rk_i^rP_{i,r}.(v'\ot v'')
&=\sum_{s=0}^r\epsilon_i^{-2s(r-s)+r^2}k_i^{r-s}(x_{i,0}^+)^{(s)}(x_{i,1}^-)^{(s)}.v'\ot k_i^s(x_0^+)^{(r-s)}(x_{i,1}^-)^{(r-s)}.v''\\
&=\sum_{s=0}^r\epsilon_i^{-2s(r-s)+r^2}(-1)^s\epsilon_i^{-s^2}k_i^{r}P_{i,s}.v'
\ot (-1)^{r-s}\epsilon_i^{-(r-s)^2}k_i^{r}P_{i,r-s}.v''.\endalign$$
Hence,
$$P_{i,r}.(v'\ot v'')=\sum_{s=0}^r P_{i,s}.v'\ot P_{i,r-s}.v''.$$
This proves relation (30) for $\calp_i^+$. The proof for $\calp_i^-$ is  similar.

The fact that $v'\ot v''$ is a highest weight vector now follows from (29).
\qed\enddemo

\demo{Proof of 7.5} In [8], we showed that
$$\Delta(x_{i,1}^-)-\Delta_i(x_{1}^-)\in\uqgh X_i^+\otimes\uqgh+
\uqgh\otimes\uqgh X_i^+,\tag31$$
where $X_i^+$ is the $\Bbb C(q)$-span of $\{x_{j,m}^+\,\mid\,m\in\Bbb Z,j\in I,j\ne i\}$. We claim that this implies that
$$\Delta(x_{i,1}^-)^r-\Delta_i(x_{1}^-)^r\in\uqgh X_i^+\otimes\uqgh+
\uqgh\otimes\uqgh X_i^+\tag32$$
for all $r\ge 1$. The proof is by induction on $r$. For the inductive step, we have
$$\align
\Delta(x_{i,1}^-)^{r+1}&=\Delta(x_{i,1}^-)^r\Delta(x_{i,1}^-)\\
&\in\left(\Delta_i(x_{1}^-)^r+\uqgh X_i^+\otimes\uqgh+\uqgh\otimes\uqgh X_i^+\right)\\
&\ \ \ \ \ \ \ \ \ \ \ \ \times\left(x_{i,1}^-\otimes 1+k_i\otimes x_{i,1}^-+\uqgh X_i^+\otimes\uqgh+\uqgh\otimes\uqgh X_i^+\right).\endalign$$
Thus, to prove (32), it is enough to show that
$$\align
\left(\uqgh X_i^+\otimes\uqgh+\uqgh\otimes\uqgh X_i^+\right)&(x_{i,1}^-\otimes 1+k_i\otimes x_{i,1}^-)\\
&\subseteq \uqgh X_i^+\otimes\uqgh+\uqgh\otimes\uqgh X_i^+.
\endalign$$
But this is clear, since by the relations in 1.2,
$$X_i^+x_{i,1}^-\subseteq \uqgh X_i^+,\ \ X_i^+k_i\subseteq \uqgh X_i^+.$$

Now, since $\uqresh$ is a $\Bbb C[q,q^{-1}]$-Hopf subalgebra of $\uqgh$, it is clear that
$$\Delta((x_{i,1}^-)^{(r)})-\Delta_i((x_{1}^-)^{(r)})\in\uqresh\otimes\uqresh.$$
On the other hand, by 6.1,
$$\uqresh=\uqreshm\uqresho+\uqresh X^{(+)}.$$
Thus, by (32), to prove the lemma, it is enough to show that
$$\left(\uqreshm\uqresho\otimes\uqreshm\uqresho\right)\cap
\left(\uqgh X^+\otimes \uqgh+\uqgh\otimes\uqgh X^+\right)=0,\tag33$$
where $X^+$ is the $\Bbb C(q)$-span of $\{x_{j,m}^+\,\mid\, j\in I,m\in\Bbb Z\}$. But (33) is a straightforward consequence of the Poincar\'e--Birkhoff--Witt basis of $\uqgh$ given in [1], Proposition 6.1.\qed\enddemo

\vskip6pt The following lemma will be needed in Section 8.

\proclaim{Lemma 7.6} Let $i\in I$ and let $\Lambda:\uepresho\to\Bbb C$ be an algebra homomorphism such that $V(\Lambda)$ is finite-dimensional. Then, $\varphi_i(U_{\epsilon_i}^{\roman{res}}(\hat{sl}_2)).v_\Lambda\subset V(\Lambda)$ is an irreducible representation of $U_{\epsilon_i}^{\roman{res}}(\hat{sl}_2)$ isomorphic to $V(\Lambda_i)$, where $\Lambda_i=\Lambda\circ\varphi_i$. \endproclaim
\demo{Proof} Let $W$ be a proper irreducible $U_{\epsilon_i}^{\roman{res}}(\hat{sl}_2)$-subrepresentation of $\varphi_i(U_{\epsilon_i}(\hat{sl}_2)).v_\Lambda$. Let $\lambda\in P^+$ be such that $\Lambda(k_i)=\epsilon_i^{\lambda(\check{\alpha}_i)}$. Let $w\in W$ be any non-zero vector of maximal weight for $U_{\epsilon_i}^{\roman{res}}(sl_2)$ (the subalgebra of $U_{\epsilon_i}^{\roman{res}}(\hat{sl}_2)$ generated by $k^{\pm 1}$, $\left[{k;0}\atop{\ell}\right]$ and the $(x_0^\pm)^{(r)}$ for all $r\in\Bbb N$. Then, 
by (29), 
$$(x_{i,n}^+)^{(r)}.w=0\ \ \ \ \ \text{for all $n\in\Bbb Z$, $r\ge 1$}.$$
But since the weight space $W_{\lambda-\eta}$ is zero unless $\eta\in Q^+$ is a multiple of $\alpha_i$, we also have
$$(x_{j,n}^+)^{(r)}.w=0\ \ \ \ \ \text{for all $n\in\Bbb Z$, $r\ge 1$, $j\ne i$}.$$
The space of vectors in $W$ annihilated by $(x_{j,n}^+)^{(r)}$ for all $n\in\Bbb Z$, $r\ge 1$, $j\in I$ is thus non-zero, and since it is clearly preserved by $\uepresho$, it contains a $\uepresh$-highest weight vector. Since $V(\Lambda)$ is irreducible as a representation of $\uepresh$, this vector must be a multiple of $v_\Lambda$. Hence, $v_\Lambda\in W$, giving the desired contradiction.\qed\enddemo

\vskip9pt
We now give an explicit construction of some highest weight representations of $U_\epsilon^{{\roman{res}}}(\hat{sl}_2)$.
\proclaim{Proposition 7.7} For any non-zero $a\in\Bbb C(q)$, there is an  algebra homomorphism $ev_a:U_q^{{\roman{res}}}(\hat{sl}_2)\to U_q^{{\roman{res}}}({sl}_2)$ such that
$${\text{ev}}_a(e_0^{{}\pm{}})=q^{{}\pm 1}a^{{}\pm 1}e^{{}\mp{}},\ {\text{ev}}_a(k_0)=k^{-1},\ {\text{ev}}_a(e_1^{{}\pm{}})=e^{{}\pm{}},\ {\text{ev}}_a(k_1) =k.\tag34$$
Moreover, we have, for all $n\in\Bbb Z$,
$${\text{ev}}_a(x_n^+)=q^{-n}a^nk^ne^+,\ \ {\text{ev}}_a(x_n^-)=q^{-n}a^ne^-k^n.\tag35$$
\endproclaim
See [4], Proposition 5.1, for the proof. 

It is obvious that, if $a\in\Bbb C$, ${\text {ev}}_a(U_q^{{\roman{res}}}(\hat{sl}_2))\subseteq U_q^{{\roman{res}}}({sl}_2)$, and hence that ${\text {ev}}_a$ induces a homomorphism of algebras ${\text{ev}}_a:U_\epsilon^{{\roman{res}}}(\hat{sl}_2)\to U_\epsilon^{{\roman{res}}}({sl}_2)$. From (34) we deduce that
$$\aligned
{\text{ev}}_a((x_n^+)^{(r)})&=\epsilon^{-nr^2}a^{nr}k^{nr}(e^+)^{(r)},\\
{\text{ev}}_a((x_n^-)^{(r)})&=\epsilon^{-nr^2}a^{nr}(e^-)^{(r)}k^{nr}.
\endaligned
\tag36$$
If $V$ is a  representation of $U_\epsilon^{{\roman{res}}}({sl}_2)$, we denote by $V_a$ the representation of $U_\epsilon^{{\roman{res}}}(\hat{sl}_2)$ obtained by pulling back $V$ by ${\text{ev}}_a$. If $V$ is of type I, then so is $V_a$. If $\overline{V}$ is a representation of $\overline{U}(sl_2)$, we define $\overline{V}_a$ similarly.

\proclaim{Proposition 7.8} Let $a\in\Bbb C$ be non-zero, and let $n\in\Bbb N$. If $v$ is a highest weight vector in $V(n)$, the $(n+1)$-dimensional irreducible type I representation of $U_\epsilon^{{\roman{res}}}({sl}_2)$, then $V(n)_a$ is a highest weight representation of $U_\epsilon^{{\roman{res}}}(\hat{sl}_2)$ with highest weight $\Lambda$ given by
$$\align
\Lambda(k)=\epsilon^n,\ \ &\Lambda\left(\left[{k;0}\atop \ell \right]\right)=
\left[n\atop\ell\right]_\epsilon,\\
\Lambda(P_r)&=(-1)^ra^r\left[n\atop{|r|}\right]_\epsilon.\tag37
\endalign$$
In particular, $P_r.v=0$ if $|r|>n$.
\endproclaim
\demo{Proof} Let $0\ne v\in V(n)_a$ be a $U_\epsilon^{{\roman{res}}}({sl}_2)$-highest weight vector. By (29),
$$(x_m^+)^{(s)}.v=0\ \ \text{for all $s\in\Bbb N$, $m\in\Bbb Z$.}$$
To prove the proposition, it therefore suffices to prove that
$$P_r.v=(-1)^ra^r\left[n\atop{|r|}\right]_\epsilon v\ \ \ \text{for all $r\in\Bbb Z$.}$$

Assume that $r>0$. Proceeding as in the proof of 7.4,
$$\align (-1)^r\epsilon^{-r^2}k^rP_r.v&=(x_0^+)^{(r)}(x_1^-)^{(r)}.v\\&=\epsilon^{-r^2+nr}(e^+)^{(r)}(e^-)^{(r)}.v\ \ \ \ \ \text{(by (36))}\\
&=\epsilon^{-r^2+nr}a^r\left[k;0\atop r\right].v\ \ \ \ \ \text{(by (2))}\\
&=\epsilon^{-r^2+nr}a^r\left[n\atop r\right]_\epsilon v.\endalign$$
Hence, $$P_r.v=(-1)^ra^r\left[n\atop r\right]_\epsilon v.$$

Finally, applying $\Omega$ to both sides of the identity in 5.1 (working in $U_q^{{\roman{res}}}(\hat{sl}_2)$ and then specialising), we obtain
$$(-1)^r\epsilon^{r^2}k^{-r}P_{-r}.v=(x_{-1}^+)^{(r)}(x_0^-)^{(r)}.v.$$
Repeating the argument in the preceding paragraph, one finds that
$$P_{-r}.v=(-1)^ra^{-r}\left[n\atop r\right]_\epsilon v.\qed$$
\enddemo
\vskip6pt\noindent{\it Remark} From equation (37), it is easy to show that
$$\Lambda(\calp^\pm(u))=\prod_{s=1}^n(1-\epsilon^{n+1-2s}a^{\pm 1}u).$$

\vskip36pt\noindent{\bf 8. Classification} 
\vskip12pt\noindent  We begin the classification of the finite-dimensional irreducible representations of $\uepresh$ with

\proclaim{Proposition 8.1} Every finite-dimensional irreducible type I representation of $\uepresh$ is highest weight.\endproclaim
\demo{Proof} Let $V$ be a finite-dimensional irreducible type I representation of $\uepresh$. Since $\dim(V)<\infty$, there exists a common eigenvector $0\ne v\in V$ for the action of $\uepresho$, which acts by commuting operators on $V$. In particular, there exists $\lambda\in P$ such that $v\in V_\lambda$. It follows that
$$V=\bigoplus_{\mu\in P} V_{\mu},$$
since the right-hand side is non-zero and preserved by $\uepresh$. Again since $\dim(V)<\infty$, there exists a maximal $\mu\in P$, say $\mu_0$, such that $V_{\mu}\ne 0$. The action of $\uepresho$ preserves $V_{\mu_0}$, so there exists $0\ne v'\in V_{\mu_0}$ that is a common eigenvector for $\uepresho$. Then, $v'$ is a $\uepresh$-highest weight vector. \qed\enddemo

In view of this result, the following theorem completes the classification of the finite-dimensional irreducible type I representations of $\uepresh$.

\proclaim{Theorem 8.2} The irreducible representation $V(\Lambda)$ of $\uepresh$ is finite-dimensional if and only if there exists an $I$-tuple $\bold P=(P_i)_{i\in I}\in\Pi^I$ such that 
$$\align
\Lambda(k_i)=\epsilon_i^{{\roman{deg}}(P_i)},\ \ \ 
&\Lambda\left(\left[{k_i;0}\atop\ell\right]\right)=
\left[{{\roman{deg}}(P_i)\atop\ell}\right]_{\epsilon_i},\\
\Lambda(\calp_i^+(u))=\frac{P_i(u)}{P_i(0)},\ \ \ 
&\Lambda(\calp_i^-(u))=\frac{Q_i(u)}{Q_i(0)},\tag38
\endalign$$
where $Q_i(u)=u^{{\roman{deg}}(P_i)}P_i(u^{-1})$.
\endproclaim
\vskip6pt\noindent{\it Remarks} 1. We shall often abuse notation by denoting the representation $V(\Lambda)$ determined by $\bold P$ as in the theorem by $V(\bold P)$.

2. It is clear that $V(\bold P)\cong V(\bold P')$ as representations of $\uepresh$ if and only if $P_i'$ is a (non-zero) multiple of $P_i$, for all 
$i\in I$. 
\vskip6pt
\demo{Proof of 8.2} Assume first that $\dim(V(\Lambda))<\infty$. Let $v_\Lambda$ be a $\uepresh$-highest weight vector in $V(\Lambda)$. Since $V(\Lambda)$ is, in particular, a finite-dimensional type I representation of $\uepres$, there exist $n_i\in\Bbb N$ such that
$$\Lambda(k_i)=\epsilon_i^{n_i^0},\ \ \ 
\Lambda\left(\left[{k_i;0}\atop\ell\right]\right)=n_i^1,$$
where $n_i=n_i^0+\ell n_i^1$, $0\le n_i^0<\ell$. Now, there is a 
$\Bbb C(q)$-algebra homomorphism $U_{q_i}({sl}_2)\to\uqgh$ that takes $e^+$ to $x_{i,-1}^+$, $e^-$ to $x_{i,1}^-$, and $k$ to $ck_i$ (this is easily checked using 1.2). This map clearly takes $U_{q_i}^{\roman{res}}({sl}_2)$ to $\uqresh$ and so induces a $\Bbb C$-algebra homomorphism $U_{\epsilon_i}^{\roman{res}}({sl}_2)\to\uepresh$. By considering $V(\Lambda)$ as a representation of $U_{\epsilon_i}^{\roman{res}}({sl}_2)$ via this homomorphism, we deduce that
$$(x_{i,1}^-)^{(r)}.v_{\Lambda}=0\ \ \ \ \text{if $r>n_i$}$$
and that $(x_{i,1}^-)^{(n_i)}.v_\Lambda$ is a non-zero multiple of $(x_{i,0}^-)^{(n_i)}.v_\Lambda$ (for $s_i(\Lambda)=\Lambda-n_i\alpha_i$ and by (28) the dimensions of the weight spaces for the action of $U_{\epsilon_i}^{\roman{res}}({sl}_2)$ are preserved by the action of the Weyl group). This implies that
$$(x_{i,0}^+)^{(r)}(x_{i,1}^-)^{(r)}.v_\Lambda\ 
\cases =0&\text{if $r>n_i$}\\
\ne 0&\text{if $r=n_i$.}\endcases$$ 
By 5.1, it follows that
$$P_{i,r}.v_\Lambda\ \cases =0&\text{if $r>n_i$}\\
\ne 0&\text{if $r=n_i$},\endcases$$
and hence that 
$$\calp_i^+(u).v_\Lambda=P_i(u)v_\Lambda$$
for some polynomial $P_i\in\Bbb C[u]$ of degree $n_i$.

It remains to prove that
$$\Lambda(\calp_i^-(u))=\frac{Q_i(u)}{Q_i(0)}.\tag39$$
We show first that $\Lambda(\calp_i^-(u))$ is uniquely determined by $P_i(u)$. By repeating {\it verbatim} the proof of Theorem 3.4 in [4] (with $q$ replaced by $\epsilon_i$), it follows that $P_i(u)$ determines $\Lambda(\psi_{i,n}^\pm)$ for all $n\in\Bbb Z$ according to the formula
$$\Lambda(\Psi_i^\pm(u))=\epsilon_i^{\pm{{\roman{deg}}}(P_i)}
\frac{P_i(\epsilon_i^{\mp 2}u)}{P_i(u)}.$$ 
By 3.1, it suffices to prove that $P_i$ determines $\Lambda(P_{i,-n\ell})$ for all $n>0$. We can assume that $n\ell\le n_i$, otherwise $\Lambda(P_{i,-n\ell})=0$ (this follows by an argument similar to that used in the previous paragraph).

Consider the vectors $(x_{i,m}^-)^{(n\ell)}.v_\Lambda$ for $m\in\Bbb Z$. These vectors are all non-zero, for by equation (2) (and an application of $T_{\check{\lambda}_i}$), we get
$$(x_{i,-m}^+)^{(n\ell)}(x_{i,m}^-)^{(n\ell)}.v_\Lambda=
\left[{k_i;0}\atop{n\ell}\right]v_\Lambda=\left({n_i^1}\atop{n}\right)v_\Lambda\ne 0.$$
Since $\dim(V(\Lambda))<\infty$, there exists a linear relation
$$\sum_{m=M}^{M'}a_m(x_{i,m}^-)^{(n\ell)}.v_\Lambda=0,$$
for certain scalars $a_m$. We can assume that $a_M=-1$. Then, we have
$$(x_{i,M}^-)^{(n\ell)}.v_\Lambda=\sum_{m=M+1}^{M'}a_m(x_{i,m}^-)^{(n\ell)}
.v_\Lambda.
\tag40$$
Applying $(x_{i,M}^+)^{(n\ell)}$ to both sides of equation (40) for {\hbox{$m=-M, -M+1,\ldots,-M'+1$}} and using 5.1 and the remark at the end of Section 6, we get a system of $M'-M$ linear equations for the $a_m$ with coefficients of the form $\Lambda(P_{i,r})$ for $r\ge 0$. It follows that the $a_m$ are given by certain rational functions of the coefficients of $P_i$. Applying $(x_{i,-M-1}^+)^{(n\ell)}$ to both sides of (40) now shows that $\Lambda(P_{i,-n\ell})$ is given by a rational function of the coefficients of $P_i$. 

By 7.5, it suffices to prove (39) when $\ung=sl_2$. Dropping the subscript $i$ and recalling that $\Lambda(\calp^-(u))$ is uniquely determined by $P(u)$, it  suffices to verify it in any one irreducible representation $V(\Lambda)$ with $\Lambda(\calp^+(u))=P(u)$. Further, since both sides of (39) are multiplicative on tensor products by 7.4, we can assume that ${\roman{deg}}(P)=1$. If $P=1-au$, where $a\in\Bbb C^\times$, then by 7.8 we can take $V(\Lambda)=V(1)_a$ and we have 
$$\Lambda(\calp^-(u))=1-a^{-1}u.$$
Since $Q(u)=uP(u^{-1})=u-a$, we get
$$\Lambda(\calp^-(u))=\frac{Q(u)}{Q(0)},$$
as required.

Turning now to the converse, we must show that, if $\bold P=(P_i)_{i\in I}\in\Pi^I$, and if $\Lambda$ is determined by $\bold P$ as in the statement of 8.2, the representation $V(\Lambda)$ is finite-dimensional. By 2.1, there exists a finite-dimensional representation $V_q(\Lambda)$ of $\uqgh$ over $\Bbb C(q)$, with highest weight vector $v_\Lambda$, such that $v_\Lambda$ has weight $\sum_{i\in I}{\roman{deg}}(P_i)\lambda_i$ for $\uqg$ and
$$\calp^+(u).v_\Lambda=\frac{P_i(u)}{P_i(0)}v_\Lambda.$$
Further, since each $P_i$ has coefficients in $\Bbb C$, it follows that
$$\uqresho.v_\Lambda=\Bbb C[q,q^{-1}]v_\Lambda.$$
By 6.1, 
$$W_q(\Lambda)=\uqreshm.v_\Lambda$$
is preserved by the action of $\uqresh$: in fact, $W_q(\Lambda)=\uqresh.v_\Lambda$. We shall prove that $W_q(\Lambda)$ is a $\Bbb C[q,q^{-1}]$-lattice in $V_q(\Lambda)$, i.e. that the natural map
$$W_q(\Lambda)\otimes_{\Bbb C[q,q^{-1}]}\Bbb C(q)\to V_q(\Lambda)$$
is a $\Bbb C(q)$-vector space isomorphism. 

For this, we note that $W_q(\Lambda)$ is the direct sum of its weight spaces for $\uqres$:
$$W_q(\Lambda)=\bigoplus_{\mu\in P}W_q(\Lambda)_\mu.$$
Choose a basis of $W_q(\Lambda)_\mu$ as a $\Bbb C[q,q^{-1}]$-module for each $\mu\in P$, and let $\Cal B$ be the union of these bases. It is clear that $\Cal B$ is linearly independent over $\Bbb C(q)$ (take a linear relation and clear denominators), and clear too that it spans $V_q(\Lambda)$ over $\Bbb C(q)$ (because $\uqresh$ spans $\uqgh$ over $\Bbb C(q)$). This proves our assertion.

Finally, define 
$$W(\Lambda)=W_q(\Lambda)\otimes_{\Bbb C[q,q^{-1}]}\Bbb C,$$
via the homomorphism $\Bbb C[q,q^{-1}]\to\Bbb C$ that takes $q$ to $\epsilon$, and denote by $w_\Lambda$ the image of $v_\Lambda\in W_q(\Lambda)$ in $W(\Lambda)$. Obviously, $W(\Lambda)=\uepresh.w_\Lambda$ and $W(\Lambda)$ is finite-dimensional (though not necessarily irreducible). Let $M$ be a maximal 
proper $\uepresh$-subrepresentation of $W(\Lambda)$, and set $V=W(\Lambda)/M$. Clearly, $w_\Lambda$ has weight $\sum_{i\in I}{\roman{deg}}(P_i)\lambda_i$ for $\uepres$, so the first two equations in the statement of 8.2 hold, and the third holds because it holds in $V_q(\Lambda)$. The last equation is proved by the argument used earlier in the proof. Hence, $V\cong V(\Lambda)$, and so $V(\Lambda)$ is finite-dimensional.\qed\enddemo

The following result is an immediate consequence of 7.4.

\proclaim{Proposition 8.3} Let $\bold P,\bold P'\in\Pi^I$, and assume that $V(\bold P)\otimes V(\bold P')$ is irreducible as a representation of $\uepresh$. Then,
$$V(\bold P)\otimes V(\bold P')\cong V(\bold P\otimes\bold P').$$
\endproclaim

\proclaim{Corollary 8.4} Let $V, V'$ be finite-dimensional irreducible type I representations of $\uepresh$, and assume that $V\otimes V'$ is irreducible. Then,
$V\otimes V'\cong V'\otimes V$ as representations of $\uepresh$.
\endproclaim
\demo{Proof} This follows from 8.1, 8.2 and 8.3.\qed\enddemo

\vskip36pt\noindent{\bf 9. A Factorization Theorem}
\vskip12pt\noindent In this section, we prove an analogue for $\uepresh$ of the factorization theorem 7.1 for finite-dimensional irreducible representations of $\uepres$. 

If $P\in\Pi$, let 
$$P=P^0P^1$$
be a factorisation such that
\vskip6pt\noindent(i) $P^0(u)$ is not divisible by $1-au^\ell$ for any non-zero $a\in\Bbb C$;

\noindent(ii) $P^1(u)=R(u^\ell)$ for some $R\in\Pi$. 
\vskip6pt\noindent Such a factorisation obviously exists and is unique up to constant multiples. If $\bold P=(P_i)_{i\in I}\in\Pi^I$, let $\bold P^0=(P_i^0)_{i\in I}$, $\bold P^1=(P_i^1)_{i\in I}$. 

\proclaim{Theorem 9.1} For any $\bold P=(P_i)_{i\in I}\in\Pi^I$, we have
$$V(\bold P)\cong V(\bold P^0)\otimes V(\bold P^1)$$
as representations of $\uepresh$.\endproclaim

To prove this theorem, we shall need the following results describing the structure of the factors $V(\bold P^0)$ and $V(\bold P^1)$. Let $\uepfinh$ be the subalgebra of $\uepresh$ generated by $\uepresho$ and the $x_{i,r}^\pm$ for all $i\in I$, $r\in\Bbb Z$.

\proclaim{Theorem 9.2} Let $\bold P=(P_i)_{i\in I}\in\Pi^I$ be such that each $P_i$ is not divisible by $1-au^\ell$ for any non-zero $a\in\Bbb C$. Then, $V(\bold P)$ is irreducible as a representation of $\uepfinh$. \endproclaim

\proclaim{Theorem 9.3} Let $\bold P=(P_i)_{i\in I}\in\Pi^I$, where $P_i(u)=R_i(u^\ell)$ for some $\bold R=(R_i)_{i\in I}\in\Pi^I$. Then,
$$V(\bold P)\cong\hat{\roman{Fr}}_\epsilon^*(\overline{V}(\bold R)),$$
the pull-back of the irreducible representation $\overline{V}(\bold R)$ of $\hat\ung$ by the Frobenius homomorphism $\hat{\roman{Fr}}_\epsilon:\uepresh\to \overline{U}(\hat\ung)$.\endproclaim

\vskip6pt\noindent{\it Remark} The representations $\overline{V}(\bold R)$ were all described explicitly in 2.5, so 9.1 and 9.3 reduce the study of $V(\bold P)$ for all $\bold P\in\Pi^I$ to the case $\bold P=\bold P^0$ to which 9.2 applies. For the latter case, we have an explicit desrciption of $V(\bold P)$ only when $\ung=sl_2$. This is given below (see 9.6).

\vskip6pt
\demo{Proof of 9.2} We first prove the theorem when $\ung=sl_2$. For this, we need 

\proclaim{Proposition 9.4} Let $r\in\Bbb N$, let $m_1,\ldots,m_r$ be positive integers $<\ell$, let $a_1,\ldots,a_r\in\Bbb C$ be non-zero, and let
$$V=V(m_1)_{a_1}\otimes\cdots\otimes V(m_r)_{a_r}.$$
Then, $V$ is irreducible as a representation of $U_\epsilon^{\roman{res}}(\hat{sl}_2)$ if and only if, for all $1\le s\ne t\le r$,
$$\frac{a_s}{a_t}\ne\epsilon^{\pm(m_s+m_t-2p)}\ \ \text{for all $0\le p<{\roman{min}}(m_s,m_t)$}.\tag41$$
Moreover, in this case, $V$ is irreducible as a representation of ${U}_\epsilon^{{\roman{fin}}}(\hat{sl}_2)$.\endproclaim
\demo{Proof} We proceed by induction on $r$, beginning with the case $r=2$ for which we change notation and consider $V(m)_a\otimes V(n)_b$. Assume that $m\le n$ (by 8.4, this is without loss of generality). The proof here is essentially the same as that in [4], Section 4.8. The crucial point is that, if $m<\ell$, the structure of $V(m)$ at an $\ell^{\roman{th}}$ root of unity is the \lq same\rq\ as its structure for generic $q$. Namely, $V(m)$ has a basis $\{v_0,v_1,\ldots,v_m\}$ and action of $U_\epsilon^{\roman{res}}(\hat{sl}_2)$ given by
$$k.v_r=\epsilon^{m-2r}v_r,\ \ e^+.v_r=[m-r+1]_\epsilon v_{r-1},\ \ 
e^-.v_r=[r+1]_\epsilon v_{r+1},$$
with $(e^\pm)^{(\ell)}$ acting as zero (and $v_{m+1}=v_{-1}=0$). As in [4], one checks that the only $U_\epsilon^{\roman{res}}({sl}_2)$-highest weight vectors in $V(m)\otimes V(n)$ are (scalar multiples of) the vectors
$$w_p=\sum_{r=0}^p(-1)^r\epsilon^{r(n-r+1)}[m-p+r]_\epsilon![n-r]_\epsilon!v_{p-r}\otimes v_r,$$
for $0\le p\le m$, and that, if $p>0$, $w_p$ is an eigenvector for $U_\epsilon^{\roman{res}}(\hat{sl}_2)^0$ and is annihilated by $x_r^+$ for all $r\in\Bbb Z$ if and only if
$$\frac{b}{a}=\epsilon^{m+n-2p+2}.\tag42$$
Note that $(x_r^+)^{(\ell)}.w_p=0$ for all $r\in\Bbb Z$, since this vector has weight $m+n-2p+2\ell$, which is $>m+n$. It follows that, if (42) holds for some $0<p\le m$, then $w_p$ generates a proper $U_\epsilon^{\roman{res}}(\hat{sl}_2)$-subrepresentation of $V(m)_a\otimes V(n)_b$. The duality argument used in [4] then shows that, if
$$\frac{b}{a}=\epsilon^{-(m+n-2p+2)}\tag43$$
for some $0<p\le m$, then $V(m)_a\otimes V(n)_b$ contains a proper $U_\epsilon^{\roman{res}}({sl}_2)$-subrepresentation containing $v_0\otimes v_0$. 

Conversely, if neither (42) nor (43) holds for any $0<p\le m$, the above argument shows that $V(m)_a\otimes V(n)_b$ is irreducible even under the subalgebra ${U}_\epsilon^{{\roman{fin}}}(\hat{sl}_2)$ of $U_\epsilon^{\roman{res}}({sl}_2)$, for neither it nor its dual representation contains a non-zero eigenvector of $U_\epsilon^{\roman{res}}(sl_2)^0$ that is annihilated by $x_r^+$ for all $r\in\Bbb Z$.

The proof for $r\ge 2$ given in [4], Section 4.8, can now be repeated {\it verbatim} (with $q$ replaced by $\epsilon$) to show that $V$ is irreducible under ${U}_\epsilon^{{\roman{fin}}}(\hat{sl}_2)$ if condition (41) holds. The converse follows from 8.4 and the $r=2$ case.\qed\enddemo

To see that the $sl_2$ case of 9.2 follows from 9.4, define an $\epsilon$-{\it segment of length} $m$, where $0<m<\ell$, to be a set of the form
$$S_m(a)=\{a\epsilon^{m-1},a\epsilon^{m-3},\ldots,a\epsilon^{-m+1}\},$$
for some non-zero $a\in\Bbb C$. Say that two $\epsilon$-segments are in {\it special position} if their union is an $\epsilon$-segment longer than both the given segments, and in {\it general position} otherwise. It is easy to see that two $\epsilon$-segments $S_{m_s}(a_s)$ and $S_{m_t}(a_t)$ are in general position if and only if (41) holds. Moreover, as in [4], Proposition 4.7, one shows that, if $P\in\Pi$ satisfies the conditions in 9.2, the set of roots of $P$, regarded as a set with multiplicities, can be written uniquely as a union of $\epsilon$-segments in general position. It follows that $P$ can be factorized uniquely as
$$P(u)=\prod_{s=1}^r\prod_{p=1}^{m_s}(1-a_s\epsilon^{m_s+1-2p}),$$
for some $r,m_1,\ldots,m_r\in\Bbb N$ and non-zero $a_1,\ldots,a_r\in\Bbb C$.
Hence, by 7.6 (and the remark which follows its proof), 8.3 and 9.4,
$$V(P)\cong V(m_1)_{a_1}\otimes\cdots\otimes V(m_r)_{a_r}.$$
By 9.4 again, $V(P)$ is irreducible under $U_\epsilon^{\roman{fin}}(\hat{sl}_2)$.
\vskip9pt
Returning now to the case of arbitrary $\ung$, we prove Theorem 9.2 in two steps:

\vskip6pt\noindent{\it Step I: If $\bold P$ is as in the statement of 9.2, then $V(\bold P)=\uepfinh.v_{\bold P}$, where $v_{\bold P}$ is a $\uepresh$-highest weight vector in $V(\bold P)$.}
\vskip6pt
For this, regard $V(\bold P)$ as a representation of $\uepres$ and let
$$V(\bold P)=\bigoplus_{\eta\in Q^+}V(\bold P)_{\lambda-\eta}$$
be its weight decomposition, where $\lambda=\sum_{i\in I}{\roman{deg}}(P_i)\lambda_i\in P^+$ is the $\uepres$-weight of $v_{\bold P}$. Note that $V(\bold P)_\lambda=\Bbb Cv_{\bold P}$. We show, by induction on $\eta$, that
$$V(\bold P)_{\lambda-\eta}\subseteq\uepfinh.v_{\bold P}.\tag44$$
We claim that this follows from the identities
$$\sum_{m=0}^n(-1)^m\epsilon_i^{-m(n+a_{ij}-1)}(x_{i,r}^-)^{(n-m)}x_{j,s}^-(x_{i,r}^-)^{(m)}=0\tag45$$
in $\uepresh$, where $n\ge 1-a_{ij}$, $r,s\in\Bbb Z$, and $i\ne j$. 

Indeed, (44) is obvious when $\eta=0$. Assuming that (44) is proved for all $\eta<\eta'\in Q^+$, we prove it for $\eta'$. By 6.1,
$$V(\bold P)_{\lambda-\eta'}=\sum_{i,r}x_{i,r}^-.V(\bold P)_{\lambda-\eta'+\alpha_i}+\sum_{i,r}(x_{i,r}^-)^{(\ell)}.V(\bold P)_{\lambda-\eta'+\ell\alpha_i}.$$
Hence, it suffices to prove that
$$(x_{i,r}^-)^{(\ell)}.V(\bold P)_{\lambda-\eta'+\ell\alpha_i}\subseteq\uepfinh.v_{\bold P}$$
for all $i\in I$, $r\in\Bbb Z$. There is nothing to prove unless $\eta'\ge\ell\alpha_i$, and if $\eta'>\ell\alpha_i$ we have
$$V(\bold P)_{\lambda-\eta'+\ell\alpha_i}=\sum_{j,s}x_{j,s}^-.V(\bold P)_{\lambda-\eta'+\ell\alpha_i+\alpha_j}$$
by the induction hypothesis again. But then, by (45),
$$\!\!\!\!\!\!\!\!\!\!\!\!(x_{i,r}^-)^{(\ell)}.V(\bold P)_{\lambda-\eta'+\ell\alpha_i}=\sum_{j,s}(x_{i,r}^-)^{(\ell)}x_{j,s}^-.V(\bold P)_{\lambda-\eta'+\ell\alpha_i+\alpha_j}\subseteq x_{i,r}^-.V(\bold P)_{\lambda-\eta'+\alpha_i}+x_{j,s}^-.V(\bold P)_{\lambda-\eta'+\alpha_j},$$
and this is contained in $\uepfinh$ by the induction hypothesis once again. 

We are thus reduced to proving that
$$(x_{i,r}^-)^{(\ell)}.v_{\bold P}\in\uepfinh.v_{\bold P}$$
for all $i$, $r$. But this is clear, since by 7.6 and the $sl_2$ case of 9.2 already proved, we have
$$(x_{i,r}^-)^{(\ell)}.v_{\bold P}\in\varphi_i(U_{\epsilon_i}^{\roman{fin}}(\hat{sl}_2)).v_{\bold P}\subseteq\uepfinh.v_{\bold P}.$$

To complete Step I, we must therefore prove the identities (45). We work in $\uqresh$ and then specialise. 

We proceed by induction on $n$. If $n=1-a_{ij}$, (45) is a special case of one of the defining relations in 1.2. Assuming that (45) holds for $n$, we multiply both sides of (45) on the left by $x_{i,r}^-$ and subtract $q_i^{-2n-a_{ij}}$ times the identity obtained by multiplying both sides of (45) on the right by $x_{i,r}^-$. The coefficient of $(x_{i,r}^-)^{(n+1-m)}x_{j,s}^-(x_{i,r}^-)^{(m)}$ on the left-hand side of the resulting identity is
$$\!\!\!\!\!\!\!\!
(-1)^m\left\{q_i^{-m(n-1+a_{ij})}[n-m+1]_{q_i}+q_i^{-2n-a_{ij}-(m-1)(n-1+a_{ij})}[m]_{q_i}\right\}=(-1)^mq_i^{-m(n+a_{ij})}[n+1]_{q_i}.$$
After cancelling the factor $[n+1]_{q_i}$, we thus have the identity (45) for $n+1$.

This completes Step I.

\vskip12pt\noindent{\it Step II: Let $\bold P$ be as in the statement of 9.2, and let $v\in V(\bold P)$ be annihilated by $x_{i,r}^+$ for all $i\in I$, $r\in\Bbb Z$. Then, $v$ is a multiple of $v_{\bold P}$.}
\vskip6pt
Note that we have already proved this when $\ung=sl_2$ (in the course of proving 9.2 itself when $\ung=sl_2$). From now on let $\ung$ be arbitrary and suppose for a contradiction that $v\in V(\bold P)$ is annihilated by $x_{i,r}^+$ for all $i\in I$, $r\in\Bbb Z$ and is not a multiple of $v_{\bold P}$. We can assume that $v$ is a weight vector for $\uepres$ and that, among vectors with these properties, $v$ has maximal weight for $\uepres$. By the analogue of (45) with $x^+$ replacing $x^-$ (and with $n=\ell$), we have
$$x_{j,s}^+(x_{i,r}^-)^{(\ell)}.v=0$$
for all $i,j\in I$, $r,s\in\Bbb Z$. The maximal property of $v$ implies that, for each $i\in I$, $r\in\Bbb Z$, $(x_{i,r}^+)^{(\ell)}.v$ is a multiple of $v_{\bold P}$ (possibly zero). We cannot have $(x_{i,r}^+)^{(\ell)}.v=0$ for all $i\in I$, $r\in\Bbb Z$, otherwise $v$ would be a multiple of $v_{\bold P}$ by the irreducibility of $V(\bold P)$. Hence, there is a unique index $i\in I$, say $i_0$, such that $(x_{i,r}^+)^{(\ell)}.v$ is a non-zero multiple of $v_{\bold P}$. It follows that $v$ has weight $\lambda-\ell\alpha_{i_0}$ for $\uepres$, where $\lambda=\sum_{i\in I}{\roman{deg}}(P_i)\lambda_i$. By 6.1, $v$ belongs to $\varphi_{i_0}(U_{\epsilon_{i_0}}^{\roman{res}}(\hat{sl}_2)).v_{\bold P}$, which is isomorphic as a representation of $U_{\epsilon_{i_0}}^{\roman{res}}(\hat{sl}_2)$ to $V(P_{i_0})$ by 7.6. By the $sl_2$-case of Step II, $v$ is a multiple of $v_{\bold P}$. This contradiction completes Step II. 
\vskip9pt
We can now complete the proof of 9.2. If $V(\bold P)$ is reducible under $\uepfinh$, it has an irreducible $\uepfinh$-subrepresentation $W$, say. Any non-zero vector $v\in W$ of maximal weight for $U_\epsilon^{\roman{res}}(\ung)$ is obviously annihilated by $x_{i,r}^+$ for all $i\in I$, $r\in\Bbb Z$. By Step II, $v$ is a multiple of $v_{\bold P}$. This contradicts Step I.\qed\enddemo

We turn now to the proof of 9.3. We shall need the following lemma.

\proclaim{Lemma 9.5} Let $i\in I$, $r\in\Bbb Z$, $n\in\Bbb N$. Then,
$$\hat{\roman{Fr}}_\epsilon((x_{i,r}^\pm)^{(n)})=\cases \frac{(\ox_{i,r}^\pm)^{n/\ell}}{(n/\ell)!}&\text{if $\ell$ divides $n$}\\
0&\text{otherwise}.\endcases$$
\endproclaim
\demo{Proof} It suffices to prove that the following diagram commutes for all $i\in\hat I$:
$$\CD
\uepresh @>\hat{\roman{Fr}}_\epsilon>>  \overline{U}(\ungh)\\
@V T_i VV     @VV \overline{T}_i V\\
\uepresh @>> \hat{\roman{Fr}}_\epsilon > \overline{U}(\ungh)
\endCD
\tag46$$
Indeed, if $T_i$ and $\overline{T}_i$ are replaced by an automorphism $\tau\in\Tau$, the corresponding diagram obviously commutes. It follows that, if (46) commutes, then
$$\hat{\roman{Fr}}_\epsilon\circ T_{\check{\lambda}_i}=\overline{T}_{\check{\lambda}_i}\circ \hat{\roman{Fr}}_\epsilon$$
for all $i\in I$. This clearly implies the lemma, in view of the isomorphism in 1.2.

To prove that the diagram (46) commutes, it suffices to prove that the homomorphisms $\hat{\roman{Fr}}_\epsilon\circ T_i$ and $\overline{T}_i\circ \hat{\roman{Fr}}_\epsilon$ agree on the generators $k_j^{\pm 1}$ and $(e_j^\pm)^{(n)}$ of $\uepresh$, for all $j\in\hat I$, $n\ge 1$. This is easy for $k_j^{\pm 1}$. For $(e_j^+)^{(n)}$, we have (see Section 1)
$$T_i((e_j^+)^{(n)})=\sum_{r=0}^{-na_{ij}}(-1)^{r-na_{ij}}\epsilon_i^{-r}
(e_i^+)^{(-na_{ij}-r)}(e_j^+)^{(n)}(e_i^+)^{(r)}\tag47$$
if $i\ne j$ (the proof when $i=j$ is similar but easier). Applying $\hat{\roman{Fr}}_\epsilon$ to the $r^{\roman{th}}$ term in the sum gives zero unless $n$ and $r$ are both divisible by $\ell$, say $r=\ell r^1$, $n=\ell n^1$. Hence, $\hat{\roman{Fr}}_\epsilon(T_i(e_j^+)^{(n)})=0$ unless $n$ is divisible by $\ell$, and 
$$\align
\hat{\roman{Fr}}_\epsilon(T_i(e_j^+)^{(\ell n^1)})&=
\sum_{r^1=0}^{-n^1a_{ij}}(-1)^{\ell r^1-\ell n^1a_{ij}}\epsilon_i^{-\ell r^1}\frac{(\oe_i^+)^{-n^1a_{ij}-r^1}}{(-n^1a_{ij}-r^1)!}\frac{(\oe^+)^{n^1}}{n^1!}\frac{(\oe_j^+)^{r^1}}{r^1!}\\
&=\sum_{r^1=0}^{-n^1a_{ij}}(-1)^{r^1-n^1a_{ij}}\frac{(\oe_i^+)^{-n^1a_{ij}-r^1}}{(-n^1a_{ij}-r^1)!}\frac{(\oe^+)^{n^1}}{n^1!}\frac{(\oe_j^+)^{r^1}}{r^1!},\endalign$$
since $\ell$ is odd and $\epsilon^\ell=1$. By the classical analogue of (47), this last sum is equal to
$$\overline{T}_i\left(\frac{(\oe_j^+)^{n^1}}{n^1!}\right)
=\overline{T}_i\circ\hat{\roman{Fr}}_\epsilon((e_j^+)^{(\ell n^1)}).$$
The proof for $(e_j^-)^{(n)}$ is similar.\qed

We can now give the

\demo{Proof of 9.3} We note first that it suffices to prove the theorem when $\ung=sl_2$. Indeed, this follows from 7.6 and its classical analogue, and the fact that
$$\hat{\roman{Fr}}_\epsilon(\varphi_i(U_{\epsilon_i}^{\roman{res}}(\hat{sl}_2)))=\overline{\varphi}_i(\overline{U}(\hat{sl}_2)),$$
which follows from 9.5.

We take $\ung=sl_2$ in the remainder of the proof. By 2.5, if $R\in\Pi$ has factorisation
$$R(u)=\prod_t(1-a_tu)^{r_t},$$
where the $a_t$ are distinct and the $r_t$ are $\ge 1$, we have
$$\overline{V}(R)\cong \bigotimes_t\overline{V}(r_t)_{a_t}.$$
Since $\hat{\roman{Fr}}_\epsilon$ is a homomorphism of Hopf algebras, the pull-back
$$\hat{\roman{Fr}}_\epsilon^*(\overline{V}(R))\cong\bigotimes_t\hat{\roman{Fr}}_\epsilon^*(\overline{V}(r_t)_{a_t}).\tag48$$

We claim next that the diagram
$$\CD
\uepresh @>\hat{{\roman{Fr}}}_\epsilon>> {\overline U}(\hat\ung)\\
@V{\roman{ev}}_bVV   @VV\overline{{\roman{ev}}}_{b^\ell}V  \\
\uepres  @>>{\roman{Fr}}_\epsilon>   \overline{U}(\ung)
\endCD
\tag49$$
commutes. To prove this, it suffices to check that ${{\roman{ev}}}_{b^\ell}\circ \hat{{\roman{Fr}}}_\epsilon$ and ${\roman{Fr}}_\epsilon\circ {\roman{ev}}_b$ agree on the generators $(e_i^\pm)^{(r)}$ and $k_i$ of $\uepresh$, for all $r\in\Bbb N$, $i=0,1$. This is clear except for $(e_0^\pm)^{(r)}$. We have
$$\align
{{\roman{ev}}}_{b^\ell}(\hat{{\roman{Fr}}}_\epsilon((e_0^\pm)^{(r)}))
&=
\cases \frac{{{\roman{ev}}}_{b^\ell}((\oe_0^\pm)^{r/\ell})}{(r/\ell)!}\ &\text{if $\ell$ divides $r$}\\
0\ &\text{otherwise}\endcases
\\
&=
\cases b^{\pm r}\frac{(\overline{e}^\mp)^{r/\ell}}{(r/\ell)!}\ &\text{if $\ell$ divides $r$}\\
0\ &\text{otherwise,}\endcases
\endalign$$
whereas 
$$\align
{\roman{Fr}}_\epsilon({\roman{ev}}_b((e_0^\pm)^{(r)}))&=
{\roman{Fr}}_\epsilon(\epsilon^{\pm r}b^{\pm r}({e}^\mp)^{(r)})\\
&=
\cases 
\epsilon^{\pm r}b^{\pm r}\frac{(\overline{e}^\mp)^{r/\ell}}{(r/\ell)!}
\ \ &\text{if $\ell$ divides $r$,}\\
0\ \ &\text{otherwise}.
\endcases
\endalign$$
Since $\epsilon^{\pm r}=1$ if $\ell$ divides $r$, the commutativity of the diagram is proved.

It follows from (48), the commutativity of the diagram (49), and 7.1(b) that
$$\hat{\roman{Fr}}_\epsilon^*(\overline{V}(R))\cong\bigotimes_t V(\ell r_t)_{a_t^{1/\ell}},$$
where $a_t^{1/\ell}$ is any $\ell^{\roman{th}}$ root of $a_t$ (the argument shows that the isomorphism class of $V(\ell r_t)_{a_t^{1/\ell}}$ is independent of which $\ell^{\roman{th}}$ root is chosen). By 7.8, the polynomial associated to $V(\ell r_t)_{a_t^{1/\ell}}$ is
$$\prod_{s=1}^{\ell r_t}(1-\epsilon^{\ell r_t+1-2s}a_t^{1/\ell}u)=\prod_{s=1}^\ell(1-\epsilon^{1-2s}a_t^{1/\ell}u)^{r_t}=(1-a_tu^\ell)^{r_t}.$$
By 7.4, $\hat{\roman{Fr}}_\epsilon^*(\overline{V}(R))\cong V(P)$, where $P(u)=R(u^\ell)$.\qed
\enddemo

\demo{Proof of 9.1} We show first that, for all $i\in I$, $n\in\Bbb Z$,
\vskip6pt\noindent(a) $x_{i,n}^\pm$ acts on $V(\bold P^0)\ot V(\bold P^1)$ as $x_{i,n}^\pm\otimes 1$;

\noindent(b) $(x_{i,n}^\pm)^{(\ell)}$ acts on $V(\bold P^0)\ot V(\bold P^1)$ as $(x_{i,n}^\pm)^{(\ell)}\otimes 1+1\otimes (x_{i,n}^\pm)^{(\ell)}$. 
\vskip6pt\noindent For this, we recall the information about the comultiplication $\Delta$ of $\uqgh$ contained in [1], Proposition 5.3.

For $i\in I$, $r\in\Bbb N$, define
$$\lambda_{i,r}=(-1)^rq_i^{-\frac12 r(r-1)}(q_i-q_i^{-1})^r[r]_{q_i}!,
 \ \ \ {\lambda}_{i,r}^*=q_i^{\frac12 r(r-1)}(q_i-q_i^{-1})^r[r]_{q_i}!,$$
and
$$R_i=\sum_{r=0}^\infty \lambda_{i,r} T_i(e_i^+)^{(r)}\otimes T_i(e_i^-)^{(r)},\ \ 
{R}_i^*=\sum_{r=0}^\infty{\lambda}_{i,r}^* (e_i^+)^{(r)}\otimes (e_i^-)^{(r)}.$$
(We have taken into account the fact that our comultiplication in $\uqgh$ is opposite to that used in [1].) Note that, since each term in the above sums belongs to $\uqresh$, $R_i$ and ${R}_i^*$ can be specialised. If $w\in\hat W$, choose a reduced decomposition
$$w=\tau s_{i_1}\ldots s_{i_N},$$
and set
$$\align
R_w&=\tau^{\otimes 2}((T_{i_1}T_{i_2}\ldots T_{i_{N-1}})^{\otimes 2}(R_{i_N})\ldots T_{i_1}^{\otimes 2}(R_{i_2})R_{i_1}),\tag50\\
{R}_{w}^*&=(T_{i_N}^{-1}T_{i_{N-1}}^{-1}\ldots T_{i_{2}}^{-1})^{\otimes 2}({R}_{i_1}^*)\ldots (T_{i_N}^{-1})^{\otimes 2}({R}_{i_{N-1}}^*){R}_{i_N}^*.\tag51\endalign$$
Then, we have [1], for all $i\in I$, $n\ge 0$,
$$\align
\Delta(x_{i,n}^-)&=R_{n\check{\lambda}_i}^{-1}(x_{i,n}^-\otimes 1+c^nk_i^{-1}\otimes x_{i,n}^-)R_{n\check{\lambda}_i},\tag52\\
\Delta(x_{i,-n}^-)&=({R}_{n\check{\lambda}_i}^*)^{-1}(x_{i,-n}^-\otimes 1+c^{-n}k_i^{-1}\otimes x_{i,-n}^-){R}_{n\check{\lambda}_i}^*,\tag53\endalign$$
and the formulas for $\Delta(x_{i,\pm n}^+)$ can be obtained from (52) and (53) by using
$$\Omega(x_{i,\mp n}^-)=x_{i,\pm n}^+,\ \ \ \Delta\circ\Omega=\Omega^{\otimes 2}\circ\Delta^{\roman{op}},$$
where $\Delta^{\roman{op}}$ is the opposite comultiplication of $\uepresh$.

We claim that $R_{n\check{\lambda}_i}$ and ${R}_{n\check{\lambda}_i}^*$ act as the identity on $V(\bold P^0)\ot V(\bold P^1)$ for all $n\in\Bbb Z$. This clearly implies (a) and (b). Indeed, using 9.3 and setting $\ell\lambda^1=\sum_{i\in I}{\roman{deg}}(P_i^1)\lambda_i$, the weight space $V(\bold P^1)_{\ell\lambda^1-\eta}$ is non-zero only if $\eta\in Q^+$ is divisible by $\ell$, so $(x_{i,n}^-)^{(r)}$ acts as zero on $V(\bold P^1)$ unless $r$ is divisible by $\ell$. We deduce from (52) and (53) that $x_{i,n}^-$ acts on $V(\bold P^0)\ot V(\bold P^1)$ as $x_{i,n}^-\otimes 1$, and from the analogue of (11) for $x^-$ that $(x_{i,n}^-\otimes 1+k_i^{-1}\otimes x_{i,n}^-)^{(\ell)}$ acts as $(x_{i,n}^-)^{(\ell)}\otimes 1+1\otimes (x_{i,n}^-)^{(\ell)}$. The case of $x_{i,n}^+$ is similar. 

As for the claim, it is easy to see from (50) that $R_{n\check{\lambda}_i}$ is a finite product of expressions of the form
$$\tau^{\otimes 2}\sum_{r=0}^\infty\lambda_{i,r}(T_{i_1}T_{i_2}\ldots T_{i_{M-1}})(e_{i_M}^+)^{(r)}\otimes 
(T_{i_1}T_{i_2}\ldots T_{i_{M-1}})(e_{i_M}^-)^{(r)}.\tag54$$
Now, $\tau s_{i_1}s_{i_2}\ldots s_{i_{M-1}}(\alpha_{i_M})$ is a real root of $\hat\ung$, so for any $\eta\in\ell Q^+$,
$$\tau^{\otimes 2}(T_{i_1}T_{i_2}\ldots T_{i_{M-1}})(e_{i_M}^-)^{(r)}(V(\bold P^1)_{\ell\lambda^1-\eta})\subseteq V(\bold P^1)_{\ell\lambda^1-\eta+r\alpha},$$
for some root $\alpha\in Q$. By 9.3, the right-hand side is zero unless $r$ is divisible by $\ell$. Hence, all the terms in the sum (54) in which $r$ is not divisible by $\ell$ are zero. But clearly $\lambda_{i,r}=0$ if $r$ is strictly positive and divisible by $\ell$, so only the $r=0$ term in (54) survives. The argument for ${R}_{n\check{\lambda}_i}^*$ is similar.

With (a) and (b) proved, we can now show that, if $v_{\bold P^0}\in V(\bold P^0)$ and $v_{\bold P^1}\in V(\bold P^1)$ are $\uepresh$-highest weight vectors, we have
$$V(\bold P)=\uepresh.(v_{\bold P^0}\otimes v_{\bold P^1}).$$
Since $V(\bold P^0)$ is irreducible under ${U}_\epsilon^{\roman{fin}}(\hat\ung)$ by 9.3, any $w^0\in V(\bold P^0)$ is a linear combination of vectors of the form
$$x_{i_1,m_1}^-x_{i_2,m_2}^-\ldots x_{i_t,m_t}^-.v_{\bold P^0}$$
with $t\in\Bbb N$, $i_1,i_2,\ldots,i_t\in I$, $m_1,m_2,\ldots,m_t\in\Bbb Z$. By (a),
$$x_{i_1,m_1}^-x_{i_2,m_2}^-\ldots x_{i_t,m_t}^-.(v_{\bold P^0}\otimes v_{\bold P^1})=(x_{i_1,m_1}^-x_{i_2,m_2}^-\ldots x_{i_t,m_t}^-.v_{\bold P^0})\otimes v_{\bold P^1},$$
hence $w^0\otimes v_{\bold P^1}\in\uepresh.(v_{\bold P^0}\otimes v_{\bold P^1})$. It follows that
$$\uepresh.(v_{\bold P^0}\otimes v_{\bold P^1})\supseteq V(\bold P^0)\otimes v_{\bold P^1}.$$
On the other hand, by 6.1 and the fact that $x_{i,n}^-$ acts as zero on $V(\bold P^1)$ for all $i\in I$, $n\in\Bbb Z$, any vector $w^1\in V(\bold P^1)$ is a linear combination of vectors of the form
$$(x_{j_1,n_1}^-)^{(\ell)}(x_{j_2,n_2}^-)^{(\ell)}\ldots (x_{j_u,n_u}^-)^{(\ell)}.v_{\bold P^1},$$
with $u\in\Bbb N$, $j_1,j_2,\ldots,j_u\in I$, $n_1,n_2,\ldots,n_u\in\Bbb Z$. We prove by induction on $u$ that
$$\uepresh.(v_{\bold P^0}\otimes v_{\bold P^1})\supseteq V(\bold P^0)\otimes (x_{j_1,n_1}^-)^{(\ell)}(x_{j_2,n_2}^-)^{(\ell)}\ldots (x_{j_u,n_u}^-)^{(\ell)}.v_{\bold P^1}.$$
If $u=0$, there is nothing to prove. Assume the result for $u-1$, and let $z^0\in V(\bold P^0)$. By (b),
$$\align
z^0\otimes (x_{j_1,n_1}^-)^{(\ell)}&(x_{j_2,n_2}^-)^{(\ell)}\ldots (x_{j_u,n_u}^-)^{(\ell)}.v_{\bold P^1}\\
&=(x_{j_1,n_1}^-)^{(\ell)}.(z^0\otimes (x_{j_2,n_2}^-)^{(\ell)}\ldots (x_{j_u,n_u}^-)^{(\ell)}.v_{\bold P^1})\\
&\ \ \ \ \ \ \ \ \ -((x_{j_1,n_1}^-)^{(\ell)}.z^0)\otimes(x_{j_2,n_2}^-)^{(\ell)}\ldots (x_{j_u,n_u}^-)^{(\ell)}.v_{\bold P^1}.
\endalign$$
By the induction hypothesis, both terms on the right-hand side of this equation belong to $\uepresh.(v_{\bold P^0}\otimes v_{\bold P^1})$, hence so does the left-hand side, completing the inductive step. 

We have now proved that $V(\bold P)$ is generated by $v_{\bold P^0}\otimes v_{\bold P^1}$ as a representation of $\uepresh$. If $V(\bold P)$ is reducible, it contains a proper irreducible subrepresentation $W$, say, and hence, by 8.1, a $\uepresh$-highest weight vector $v\ne 0$. We can write
$$v=\sum_p v_p^0\otimes v_p^1,$$
where $v_p^0\in V(\bold P^0)$, $v_p^1\in V(\bold P^1)$ and the $v_p^1$ are linearly independent. By (a), for all $n\in\Bbb Z$,
$$0=x_n^+.v=\sum_p (x_n^+.v_p^0)\otimes v_p^1,$$
so $x_n^+.v_p^0=0$ for all $p$. Since $V(\bold P^0)$ is irreducible as a representation of ${U}_\epsilon^{\roman{fin}}(\hat\ung)$ by 9.3, this implies that each $v_p^0$ is a scalar multiple of $v_{\bold P^0}$. (For, the space of vectors in $V(\bold P^0)$ annihilated by $x_{i,n}^+$ for all $i\in I$, $n\in\Bbb Z$ is preserved by $\uepresho$, hence contains a $\uepresho$-eigenvector.) So
$$v=v_{\bold P^0}\otimes z^1$$
for some $z^1\in V(\bold P^1)$. By (b), for all $i\in I$, $n\in\Bbb Z$,
$$0=(x_{i,n}^+)^{(\ell)}.v=v_{\bold P^0}\otimes (x_{i,n}^+)^{(\ell)}.z^1,$$
so $(x_{i,n}^+)^{(\ell)}.z^1=0$. Since $x_{i,n}^+$ acts as zero on $V(\bold P^1)$ for all $i\in I$, $n\in\Bbb Z$, this forces $z^1$ to be a multiple of $v_{\bold P^1}$. But then $v_{\bold P^0}\otimes v_{\bold P^1}\in W$. In view of the first part of the proof, this contradicts the fact that $W$ is a proper subrepresentation.\qed\enddemo

We conclude by giving an explicit realization of all the finite-dimensional irreducible type I representations of $U_\epsilon^{\roman{res}}(\hat{sl}_2)$:

\proclaim{Theorem 9.6} Every finite-dimensional irreducible type I representation of $U_\epsilon^{\roman{res}}(\hat{sl}_2)$ is isomorphic to a tensor product 
$$V(m_1)_{a_1}\otimes\cdots\otimes V(m_r)_{a_r}\otimes V(\ell n_1)_{b_1}\otimes \cdots\otimes V(\ell n_s)_{b_s},\tag55$$
where $r,s,0\le m_1,\ldots,m_r<\ell,n_1,\ldots,n_s\in\Bbb N$ and $a_1,\ldots,a_r,b_1,\ldots,b_s\in\Bbb C$ are non-zero. Two such irreducible tensor products, with parameters $r,s,m_1,\ldots,m_r$, $n_1,\ldots,n_s,a_1,\ldots,a_r,b_1,\ldots,b_s$ and $r',s',m_1',\ldots,m_{r'}',n_1',\ldots,n_{s'}',a_1',\ldots,a_{r'}'$, 

\noindent $b_1',\ldots,b_{s'}'$, respectively, are isomorphic if and only if
\vskip6pt\noindent (i) $r=r'$ and $s=s'$ and

\noindent(ii) there are permutations $\rho$ of $\{1,2,\ldots,r\}$ and $\sigma$ of $\{1,2,\ldots,s\}$ such that
$$m_t'=m_{\rho(t)},\ \ a_t'=a_{\rho(t)},\ \ n_t'=n_{\sigma(t)},\ \ (b_t')^\ell=b_{\sigma(t)}^\ell\ \ \text{for all $t$}.$$

Conversely, a tensor product (55), where $r,s,0\le m_1,\ldots,m_r<\ell$, 

\noindent $n_1,\ldots,n_s\in\Bbb N$ and $a_1,\ldots,a_r,b_1,\ldots,b_s\in\Bbb C$ are non-zero, is irreducible if and only if
\vskip6pt\noindent (i) for all $1\le t\ne u\le r$,
$$\frac{a_t}{a_u}\ne\epsilon^{\pm(m_t+m_u-2p)}\ \ \text{for all $0\le p<{\roman{min}}(m_t,m_u)$, and}$$
(ii) $b_1^\ell,\ldots,b_s^\ell$ are distinct.
\endproclaim
\demo{Proof} This is a straightforward consequence of 8.2 and 9.1--9.4.\qed\enddemo

\vskip36pt\noindent{\bf References} 
\vskip12pt\noindent
1. Beck, J., Braid group action and quantum affine algebras, {\it Commun. Math. Phys.} {\bf 165} (1994), 555--568.

\noindent
2. Beck, J. and Kac, V. G., Finite-dimensional representations of quantum affine algebras at roots of unity, {\it J. Amer. Math. Soc.} {\bf 9} (1996), 391--423.

\noindent
3. Chari, V. and Pressley, A. N., New unitary representations of loop groups, {\it Math. Ann.} {\bf 275} (1986), 87--104.

\noindent
4.  Chari, V. and Pressley, A. N.,  Quantum affine algebras, {\it Commun. Math. Phys.} {\bf 142} (1991), 261--283.

\noindent
5. Chari, V. and Pressley, A. N., {\it A Guide to Quantum Groups}, Cambridge University Press, Cambridge, 1994.

\noindent 
6.  Chari, V. and Pressley, A. N., Quantum affine algebras and their representations, Canadian Math. Soc. Conf. Proc. {\bf 16} (1995), 59--78.

\noindent 
7.  Chari, V. and Pressley, A. N., Yangians, integrable quantum systems and Dorey's rule, {\it Commun. Math. Phys.} {\bf 181} (1996), 265--302.

\noindent 
8. Chari, V. and Pressley, A. N., Minimal affinizations of representations of quantum groups: the simply-laced case, {\it J. Algebra} {\bf 184} (1996), 1--30.

\noindent 
9. Chari, V. and Pressley, A. N., Quantum affine algebras and rationality, to appear in {\it Proceedings of the NATO Advanced Study Institute, Cargese, 1996}, Plenum Press, 1997.

\noindent
10. De Concini, C. and Kac, V. G., Representations of quantum groups at roots of 1, in {\it Operator Algebras, Unitary Representations, Enveloping Algebras and Invariant Theory}, A. Connes, M. Duflo, A. Joseph and R. Rentschler (eds.), pp. 471--506, Progress in Mathematics {\bf 92}, Birkh\"auser, Boston, 1990. 

\noindent
11. De Concini, C., Kac, V. G. and Procesi, C., Quantum coadjoint action, {\it J. Amer. Math. Soc.} {\bf 5} (1992), 151--190.

\noindent
12. Drinfel'd, V. G., A new realization of Yangians and quantized affine algebras, {\it Soviet Math. Dokl.} {\bf 36}, 212--216.

\noindent
13. Garland, H., The arithmetic theory of loop algebras, {\it J. Algebra} {\bf 53} (1978), 480--551.

\noindent 
14. Jing, N.-H., On Drinfeld realization of quantum affine algebras, preprint q-alg/9610035.

\noindent
15. Lusztig, G. {\it Introduction to Quantum Groups}, Progress in Mathematics {\bf 110}, Birkh\"auser, Boston, 1993.

\enddemo

\vskip36pt
{\eightpoint{
$$\matrix\format\l&\l&\l&\l\\
\phantom{.} & {\text{Vyjayanthi Chari}}\phantom{xxxxxxxxxxxxx} & {\text{Andrew\ Pressley}}\\
\phantom{.}&{\text{Department of Mathematics}}\phantom{xxxxxxxxxxxxx} & {\text{Department of Mathematics}}\\
\phantom{.}&{\text{University of California}}\phantom{xxxxxxxxxxxxx} & {\text{King's College}}\\
\phantom{.}&{\text{Riverside}}\phantom{xxxxxxxxxxxxx} & {\text{Strand}}\\
\phantom{.}&{\text{CA 92521}}\phantom{xxxxxxxxxxxxx} & {\text{London WC2R 2LS}}\\
\phantom{.}&{\text{USA}}\phantom{xxxxxxxxxxxxx} & {\roman{UK}}\\
&{\text{email: chari\@ucrmath.ucr.edu}}\phantom{xxxxxxxxxxxxx} & {\text{email:anp\@uk.ac.kcl.mth}}
\endmatrix$$
}}

\enddocument